\documentclass[twoside, 10pt]{article}

\usepackage{amsmath}
\usepackage{paralist}
\usepackage{amssymb,mathrsfs,psfrag,graphicx}
\usepackage{hyperref}
\def\v{\varepsilon}
\def\e{\epsilon}

\def\t{\theta}
\def\T{\Theta}
\def\k{\kappa}

\def\a{\alpha}

\def\g{\gamma}
\def\d{\delta}

\def\s{\sigma}
\def\z{\zeta}

\def\di{\displaystyle}

\newtheorem{theorem}{Theorem}[section]

\newtheorem{lemma}[theorem]{Lemma}
\newtheorem{proposition}{Proposition}[section]

\newtheorem{remark}{Remark}

\begin{document}

\title{\bf   Stability of Superposition of Viscous Contact Wave and Rarefaction Waves
for Compressible Navier-Stokes System} \vskip 0.5cm
\author{
 \quad Feimin Huang\thanks{Academy of Mathematics and Systems
   Science, Chinese Academy of Sciences, Beijing 100190, P. R. China and
  Beijing Center of Mathematics
and Information Sciences, Beijing 100048, P. R. China
   ({\tt fhuang@amt.ac.cn}). The work of F. Huang is partially supported by the NSFC Grants No. 11371349,
   National Basic Research Program of China (973 Program) under Grant No. 2011CB808002.},
\quad Teng Wang\thanks{ Academy of Mathematics and
Systems Science, Chinese Academy of Sciences, Beijing 100190, P. R. China
   ({\tt tengwang@amss.ac.cn}). }
 }

\date{}
\maketitle
\begin{abstract}

This paper is concerned with the large-time behavior of solutions for the one-dimensional
compressible Navier-Stokes system. We show that the combination of viscous contact wave with rarefaction waves
for the non-isentropic  polytropic gas is stable under \emph{large} initial perturbation
without the condition that the adiabatic exponent $\gamma$ is close to 1, provided the strength of the combination waves is suitably small.


\


 \noindent{\it Key
words and phrases:} viscous contact discontinuity, compressible Navier-Stokes system,
 stability, large initial perturbation
\end{abstract}
\section{Introduction}
\renewcommand{\theequation}{\arabic{section}.\arabic{equation}}
\setcounter{equation}{0}

The one-dimensional compressible Navier-Stokes system in Lagrangian coordinates read
\begin{equation}\label{ns}
\left\{
\begin{array}{ll}
\di v_t-u_x=0,\\
\di u_t+p_x=\mu\left(\frac{u_x}{v}\right)_x,\\
\di \left(e+\frac{u^2}{2}\right)_t+\left(p u\right)_x=\left(\kappa \frac{\t_x}{v}+\mu\frac{uu_x}{v}\right)_x
\end{array}
\right.
\end{equation}
for $x\in\mathbb{R}=(-\infty,+\infty)$, $t>0$, where $v(x,t)>0$, $u(x,t)$, $\t(x,t)>0$, $e(x,t)>0$ and $p(x,t)$
are the specific volume, fluid velocity, absolute temperature, internal energy
and pressure, respectively, while the positive constants $\mu$ and $\k$ denote
the viscosity and heat conduction coefficients, respectively. Here we study the ideal
polytropic fluids so that $p$ and $e$ are given by the state equations
\begin{equation*}
\di p=\frac{R\t}{v}=Av^{-\gamma}\exp\left(\frac{\gamma-1}{R}s\right),\quad e=c_{\nu}\t+\mathrm{const}.,
\end{equation*}
where $s$ is the entropy, $\gamma>1$ is the adiabatic exponent,  $c_{\nu}=\frac{R}{\gamma-1}$ is
the specific heat, and $A$ and $R$ are both positive constants.
We consider the Cauchy problem to the system \eqref{ns} supplement with the following initial
and far field conditions:
\begin{equation}\label{initial}
\left\{
\begin{array}{ll}
\di (v,u,\t)(x,0)=(v_0,u_0,\t_0)(x), &\di x\in \mathbb{R},\\
\di (v,u,\t)(\pm\infty,t)=(v_{\pm},u_{\pm},\t_{\pm}), &\di t>0,
\end{array}
\right.
\end{equation}
where $v_{\pm}(>0)$, $u_{\pm}$ and $\t_{\pm}(>0)$ are given constants, and we assume
$\inf_{\mathbb{R}}v_0>0$, $\inf_{\mathbb{R}}\t_0>0$, and $(v_0,u_0,\t_0)(\pm\infty)=(v_{\pm},u_{\pm},\t_{\pm})$
as compatibility conditions. When the far field states are the same, i.e., $v_+=v_-$, $u_+=u_-$, $\t_+=\t_-$,
there has been considerable progress on the global existence of the solutions to the system \eqref{ns} since 1977,
see \cite{jiang,jiang-2,kazhikhov,kazhi-she,L-L} and the reference therein. In particular, Jiang \cite{jiang, jiang-2} first obtained some
interesting results on the large-time behavior of solutions, however the temperature is only shown to be locally bounded in space. More recently,  Li and Liang \cite{L-L} improved Jiang's results by proving the temperature is uniformly bounded.

The existence and large time behavior of solutions to the system \eqref{ns} with different end states
become much more complicated. It is noted that, if the dissipation effects are neglected, i.e.,
$\mu=\k=0$, the system \eqref{ns} is reduced to the compressible Euler equations as follows
\begin{equation}\label{euler}
\left\{
\begin{array}{ll}
\di v_t-u_x=0,\\
\di u_t+p_x=0,\\
\di \left(e+\frac{u^2}{2}\right)_t+\left(p u\right)_x=0,
\end{array}
\right.
\end{equation}
which is the most important hyperbolic system of conservation laws. It is well known that
the system \eqref{euler}
has rich wave phenomena. Indeed, it contains three basic wave patterns (see \cite{Smoller}), two nonlinear waves:
shock and rarefaction wave, and a linearly degenerate wave: contact discontinuity. When we consider the Riemann initial data
\begin{equation}\label{Riemann}
(v,u,\t)(x,0)=\left\{
\begin{array}{ll}
\di (v_-,u_-,\t_-),&\di x<0,\\
\di (v_+,u_+,\t_+),&\di x>0,
\end{array}
\right.
\end{equation}
the solutions consist of the above three wave patterns
and their superpositions, called by Riemann solutions, and govern both the local and large time
asymptotic behavior of general solutions of the system \eqref{euler}. It is of great importance and interest
to study the large-time behavior of the viscous version of these basic wave patterns and their superpositions
to the compressible Navier-Stokes system \eqref{ns}.

There has been extensive literature on the stability analysis of viscous wave pattern to system \eqref{ns},
meanwhile new phenomena has been discovered and new techniques have been developed. We refer to
\cite{Huang-matsumura,KM,MN-85} for the shock wave, \cite{MN-86,MN-92,nishi-yang-zhao} for the rarefaction wave,
\cite{huang-ma-shi,Huang-Matsumura-Xin,Huang-Xin-Yang,Huang-Yang,huang-zhao} for the viscous contact discontinuity, and the reference therein.
However, the stability of the superposition of several wave patterns is more complicated and challenging due to
the fact that the stability analysis essentially depends on the underlying properties of basic wave pattern and these frameworks are not compatible with each other. Besides,  the wave interaction between different families of wave patterns is complicated. Recently, Huang-Matsumura in \cite{Huang-matsumura} showed that the superposition of the two viscous shock profiles
for the Navier-Stokes system \eqref{ns} is asymptotically stable without the zero initial mass condition.
This result was extended in \cite{Huang-Li-matsumura} to the combination of
viscous contact discontinuity with rarefaction waves  by deriving new estimates on the heat kernel.
The time-asymptotic stability of other cases is still open!

It is noted that all results mentioned above are concerned with the small perturbation around the viscous
wave pattern. In another word, they are ``local" stability. A nature problem is, whether or not these basic wave patterns and their linear superpositions
are stable even for large perturbation. This is equivalent to study the global stability of the viscous
wave patterns to the system \eqref{ns}, which is challenging because the nonlinear terms play leading role in the large solutions, while the linearized system around wave patterns is essential for the local stability. Along this direction, Nishihara-Yang-Zhao in \cite{nishi-yang-zhao} first proved the
rarefaction waves for the system \eqref{ns} are stable with ``partially" large perturbation with the condition that the adiabatic exponent $\gamma$ is
closing enough to ``1". Precisely speaking, the amplitude of initial perturbation is reciprocal to $\gamma-1$. That is, $\gamma-1$ is smaller, the perturbation around the rarefaction wave can be larger. This result is extended by Huang-Zhao in \cite{huang-zhao} to the case of single viscous contact wave and the combination of viscous contact wave and rarefaction
waves for a free boundary value problem, and further by Hong in \cite{Hong} for the Cauchy problem.
Note that the condition on $\gamma$ plays essential role in \cite{Hong,huang-zhao,nishi-yang-zhao}, however, it is not natural in the physical setting. The main aim of this paper is to remove the non-physical condition of the adiabatic exponent $\gamma$.

\

Before stating the main results, we first recall the viscous contact wave $(V,U,\Theta)$
for the compressible Navier-Stokes system \eqref{ns} introduced in \cite{Huang-Xin-Yang}. For the Riemann
problem \eqref{euler}-\eqref{Riemann}, it is known that the contact discontinuity solution takes the form
\begin{equation}
(\tilde V,\tilde U,\tilde\Theta)(x,t)=\left\{
\begin{array}{ll}
\di (v_-,u_-,\t_-),&\di x<0,t>0,\\
\di (v_+,u_+,\t_+),&\di x>0,t>0,
\end{array}
\right.
\end{equation}
provide that
\begin{equation}\label{RH}
\di u_-=u_+,\quad p_-\triangleq\frac{R\theta_-}{v_-}=p_+\triangleq\frac{R\t_+}{v_+}.
\end{equation}
We assume that $u_-=u_+=0$ without loss of generality.
Due to the effect of heat conductivity, the contact discontinuity $(\tilde V,\tilde U,\tilde \Theta)$
is smoothed and behaves as a diffusion wave, called by "viscous contact wave". The viscous contact wave $(V,U,\Theta)$ can be constructed
as follows. Since the pressure for the profile $(V,U,\Theta)$ is expected to be
constant asymptotically, we set
$$
\frac{R\Theta}{V}=p_+,
$$
which indicates the leading part of the energy equation $\eqref{ns}_3$ is
\begin{equation}\label{Theta t}
\di c_{\nu}\Theta_t+p_+U_x=\k\left(\frac{\Theta_x}{V}\right)_x.
\end{equation}
The equation \eqref{Theta t} and $\eqref{ns}_1$ lead to a nonlinear diffusion equation,
\begin{equation}\label{1.8}
\di \Theta_t=a\left(\frac{\Theta_x}{\Theta}\right)_x,\quad
\Theta(\pm\infty,t)=\theta_{\pm},\quad a=\frac{\kappa p_+(\gamma-1)}{\gamma R^2}>0,
\end{equation}
which has a unique self-similar solution $\Theta(x,t)=\Theta(\xi)$, $\xi=\frac{x}{\sqrt{1+t}}$
due to \cite{xiao-liu}. Furthermore, $\Theta(\xi)$ is a monotone function,
increasing if $\theta_+>\theta_-$ and decreasing if $\theta_+<\theta_-$. On the other hand,
there exists some positive constant $\delta$, such that for $\d=|\t_+ -\t_-|$, $\Theta$ satisfies
\begin{equation}\label{1.9}
\di (1+t)|\Theta_{xx}|+(1+t)^{\frac{1}{2}}|\Theta_x|+|\Theta-\theta_{\pm}|
=O(1)\delta e^{-\frac{c_1x^2}{1+t}} \quad \mathrm{as}~ |x|\rightarrow\infty,
\end{equation}
where $c_1$ is positive constant depending only on $\t_{\pm}$. Once $\Theta$ is determined,
the contact wave profile $(V,U,\Theta)(x,t)$ is then defined as follows:
\begin{equation}\label{contact}
\di V=\frac{R}{p_+}\Theta,\quad U=\frac{\k(\gamma-1)}{\gamma R}\frac{\Theta_x}{\Theta},
\quad \Theta=\Theta.
\end{equation}
The contact wave $(V, U,\Theta)(x,t)$ solves the compressible Navier-Stokes system \eqref{ns}
time asymptotically, that is,
\begin{equation}\label{contact equation}
\left\{
\begin{array}{ll}
\di V_t-U_x=0,\\
\di U_t+\left(\frac{R\Theta}{V}\right)_x=\mu\left(\frac{U_x}{V}\right)_x+R_1,\\
\di c_{\nu}\Theta_t+p(V,\Theta)U_x=\left(\kappa \frac{\Theta_x}{V}\right)_x+\mu\frac{U^2_x}{V}+R_2,
\end{array}
\right.
\end{equation}
where
\begin{equation}
\di \widetilde R_1=U_t-\mu\left(\frac{U_x}{V}\right)_x,\quad \widetilde R_2=-\mu\frac{U_x^2}{V}.
\end{equation}

We first study the global stability of single viscous contact wave $(V, U,\Theta)$ for arbitrary $\gamma>1$. For this, we put the perturbation $(\phi,\psi,\zeta)(x,t)$ by
\begin{equation}
(\phi,\psi,\zeta)(x,t)=(v-V,u-U,\t-\Theta)(x,t).
\end{equation}

The precise statement of the first result is
\begin{theorem}{\bf (Viscous contact wave)}\label{theorem1}
For any given left end state $(v_-,u_-,\theta_-)$, suppose that the right end state $(v_+,u_+,\t_+)$ satisfies \eqref{RH}.
Let $(V,U,\Theta)$ be the viscous contact wave defined in \eqref{contact} with strength
$\d=|\t_+-\t_-|$. There exist a function $m(\d)$ satisfying $m(\d)\to +\infty$, as $\d\to 0$ and a small constant $\delta_0$ such that if $\d<\d_0$ and the initial data satisfies
\begin{equation}\label{data}
\left\{
\begin{array}{ll}
\di v_0(x), \t_0(x)\geq m_0^{-1},\quad m_0=:m(\d_0),\\[0.2cm]
\|(v_0(x)-V(x,0),u_0(x)-U(x,0),\t_0(x)-\Theta(x,0))\|_{H^1(\mathbb{R})}\le m_0,
\end{array}
\right.
\end{equation}
then the Cauchy problem \eqref{ns}-\eqref{initial} admits
a unique global solution $(v,u,\t)$ satisfying
$$
(v-V,u-U,\theta-\Theta)(x,t)\in C\big((0,+\infty);H^1(\mathbb{R})\big);
$$
$$
(v-V)_x(x,t)\in L^2(0,+\infty;L^2\big(\mathbb{R})\big);
$$
$$
(u-U,\t-\Theta)_x(x,t)\in L^2\big(0,+\infty;H^1(\mathbb{R})\big).
$$
Furthermore,
\begin{equation}\label{behavior-1}
\lim_{t\rightarrow+\infty}\sup_{x\in\mathbb{R}}\left|(v-V,u-U,\t-\Theta)(x,t)\right|=0.
\end{equation}

\end{theorem}

\begin{remark}
Theorem 1.1 means that if the strength of contact wave is smaller, the initial perturbation can be larger. In particular, when $\d=0$, that is, the asymptotic state is a constant one $(\bar{v},\bar{u},\bar{\t})$ instead of a wave pattern, then for any initial data $(v_0-\bar{v},u_0-\bar{u},\t_0-\bar{\t})(x)\in H^1(\mathbb{R})$, there always exists a small constant $\delta_0$ such that
\eqref{data} holds, in which the second term is replaced by $\|(v_0(x)-\bar{v},u_0(x)-\bar{u},\t_0(x)-\bar{\t})\|_{H^1(\mathbb{R})}\le m_0$. This coincides with the one in Li-Liang \cite{L-L}.
\end{remark}

\begin{remark}
Theorem 1.1 holds for any $\gamma>1$ and thus removes the condition that $\gamma$ is close to 1 in
Nishihara-Yang-Zhao \cite{nishi-yang-zhao} and  also in \cite{Hong,huang-zhao}.
\end{remark}


When the relation \eqref{RH} fails, the basic theory of hyperbolic systems of conservation laws
implies that for any given constant state $(v_-, u_-, \t_-)$ with $v_->0$, $\t_->0$ and $u_-\in\mathbb{R}$,
there exists a suitable neighborhood $\Omega(v_-,u_-,\t_-)$ of $(v_-, u_-, \t_-)$ such that for any
$(v_+, u_+, \t_+)\in\Omega(v_-,u_-,\t_-)$, the Riemann problem of the Euler system \eqref{euler}, \eqref{Riemann}
has a unique solution. In this paper, we only consider the case of the superposition of the viscous
contact wave and rarefaction waves with
\begin{equation}\label{zheng}
\di (v_+,u_+,\t_+)\in R_1CR_3(v_-,u_-,\t_-)\subset\Omega(v_-,u_-,\t_-),
\end{equation}
where
\begin{equation*}
\begin{array}{ll}
\di R_1CR_3(v_-,u_-,\t_-)\triangleq\Bigg\{(v,u,\t)\in\Omega(v_-,u_-,\t_-)\Bigg|s\neq s_-,\\
\di u\geq u_--\int_{v_-}^{e^{\frac{\g-1}{R\g}(s_--s)}v}\lambda_-(\eta,s_-)d\eta,
u\geq u_--\int^v_{e^{\frac{\g-1}{R\g}(s_--s)}v_-}\lambda_+(\eta,s)d\eta\Bigg\}
\end{array}
\end{equation*}
and
$$
s=\frac{R}{\g-1}\ln\frac{R\t}{A}+R\ln v,\quad s_{\pm}=\frac{R}{\g-1}\ln\frac{R\t_{\pm}}{A}+R\ln v_{\pm},~~
\lambda_{\pm}(v,s)=\pm\sqrt{A\g v^{-\g-1}e^{\frac{\g-1}{R}s}}.
$$

By the standard argument (e.g. \cite{Smoller}), there exists a unique pair of points $(v_-^m, u^m, \t_-^m)$
and $(v_+^m, u^m, \t_+^m)$ in $\Omega(v_-,u_-,\t_-)$ satisfying
$$
\frac{R\t_-^m}{v_-^m}=\frac{R\t_+^m}{v_+^m}\triangleq p^m,
$$
the points $(v_-^m,u^m,\t_-^m)$ and $(v_+^m,u^m,\t_+^m)$ belong to the 1-rarefaction wave
curve $R_-(v_-,u_-,\t_-)$ and the 3-rarefaction wave curve $R_+(v_+,u_+,\t_+)$, respectively, where
$$
R_{\pm}(v_{\pm},u_{\pm},\t_{\pm})=\left\{(v,u,\t)\Bigg|s=s_{\pm},
u=u_{\pm}-\int_{v_{\pm}}^v\lambda_{\pm}(\eta,s_{\pm})d\eta,v>v_{\pm}
\right\}.
$$

Without loss of generality, we assume $u^m=0$ in what follows.
The 1-rarefaction wave $(v_-^r,u_-^r,\t_-^r)(\frac{x}{t})$ (respectively the 3-rarefaction wave $(v_+^r,u_+^r,\t_+^r)(\frac{x}{t})$)
connecting $(v_-,u_-,\t_-)$ and $(v_-^m,0,\t_-^m)$ (respectively $(v_+^m,0,\t_+^m)$ and $(v_+,u_+,\t_+)$) is the
weak solution of the Riemann problem of the Euler system \eqref{euler} with the following initial Riemann data
\begin{equation}\label{Riemann2}
\di(v_{\pm},u_{\pm},\t_{\pm})(x,0)=\left\{
\begin{array}{ll}
(v_{\pm}^m,0,\t_{\pm}^m), &\quad\pm x<0,\\
(v_{\pm},u_{\pm},\t_{\pm}), &\quad\pm x>0.
\end{array}
\right.
\end{equation}
Since the rarefaction wave $(v_{\pm}^r,u_{\pm}^r,\t_{\pm}^r)$ are weak solutions, it is
convenient to construct approximate rarefaction wave which is smooth. Motivated by \cite{MN-86}, the smooth solutions of Euler
system \eqref{euler}, $(V_{\pm}^r,U_{\pm}^r,\T_{\pm}^r)$, which approximate $(v_{\pm}^r,u_{\pm}^r,\t_{\pm}^r)$,
are given by
\begin{equation}\label{appro-rare}
\di\left\{
\begin{array}{ll}
\di\lambda_{\pm}(V_{\pm}^r(x,t),s_{\pm})=w_{\pm}(x,t),\\
\di U_{\pm}^r=u_{\pm}-\int_{v_{\pm}}^{V_{\pm}^r(x,t)}\lambda_{\pm}(\eta,s_{\pm})d\eta,\\
\T_{\pm}^r=\t_{\pm}(v_{\pm})^{\g-1}(V_{\pm}^r)^{1-\g},
\end{array}
\right.
\end{equation}
 where $w_-$ (respectively $w_+$) is the solution of the initial problem for the typical Burgers equation:
\begin{equation}\label{burgers}
\di\left\{
\begin{array}{ll}
\di w_t+ww_x=0,\quad(x,t)\in\mathbb{R}\times(0,\infty),\\
\di w(x,0)=\frac{w_r+w_l}{2}+\frac{w_r-w_l}{2}\tanh x,
\end{array}
\right.
\end{equation}
with $w_l=\lambda_-(v_-,s_-)$, $w_r=\lambda_-(v_-^m,s_-)$ (respectively $w_l=\lambda_+(v_+^m,s_+)$, $w_r=\lambda_+(v_+,s_+)$).

Let $(V^{cd},U^{cd},\T^{cd})(x,t)$ be the viscous contact wave constructed in \eqref{1.8} and \eqref{contact}
with $(v_{\pm},u_{\pm},\t_{\pm})$ replaced by $(v_{\pm}^m,0,\t_{\pm}^m)$, respectively.

To describe the strengths of the viscous contact wave and rarefaction waves for later use, we set
\begin{equation*}
\begin{array}{ll}
\d^{r_1}=|v_-^m-v_-|+|0-u_-|+|\t_-^m-\t_-|, \quad \d^{cd}=|\t_+^m-\t_-^m|,\\
\d^{r_3}=|v_+^m-v_+|+|0-u_+|+|\t_+^m-\t_+|
\end{array}
\end{equation*}
and $\d=\min(\d^{r_1},\d^{cd},\d^{r_3})$. If
\begin{equation}\label{same-order}
\d^{r_1}+\d^{cd}+\d^{r_3}\leq C\d,\quad {\rm as }\quad \d^{r_1}+\d^{cd}+\d^{r_3}\rightarrow 0
\end{equation}
holds for a positive constant $C$, we call the strengths of the wave patterns ``small with the same order". In this case, we have
\begin{equation}
\d^{r_1}+\d^{cd}+\d^{r_3}\leq C|(v_+-v_-,u_+-u_-,\t_+-\t_-,)|.
\end{equation}
In what follows, we always assume \eqref{same-order}. We define
\begin{equation}\label{ansatz}
\left(
\begin{array}{l}
\di V\\
\di U \\
\di \T
\end{array}
\right)(x,t) = \left(
\begin{array}{l}
\di V^{cd}+V_-^r+V_+^r\\
\di U^{cd}+U_-^r+U_+^r\\
\di \T^{cd}+\T_-^r+\T_+^r
\end{array}\right)(x,t)-\left(
\begin{array}{l}
\di v_-^m+v_+^m\\
\di \quad~~ 0\\
\di \t_-^m+\t_+^m
\end{array}\right),
\end{equation}
and
$$
(\phi,\psi,\zeta)(x,t)=(v-V,u-U,\t-\T)(x,t).
$$

The precise statement of the second result is

\begin{theorem}{\bf (Composite waves)}\label{theorem2}
For any given left end state $(v_-,u_-,\theta_-)$, let $(V,U,\Theta)$ be defined in \eqref{ansatz} with strength
satisfying \eqref{same-order}. Then there exists a function $m(\d)$ satisfying $m(\d)\to +\infty$,
as $\d\to 0$ and a small constant $\delta_0$ , such that if $|(v_+-v_-,u_+-u_-,\t_+-\t_-)|<\d_0$ and the initial data satisfies
\begin{equation}\label{data2}
\left\{
\begin{array}{ll}
\di  v_0(x), \t_0(x)\geq m_0^{-1},\quad m_0=:m(\d_0),\\[0.2cm]
\|(v_0(x)-V(x,0),u_0(x)-U(x,0),\t_0(x)-\Theta(x,0))\|_{H^1(\mathbb{R})}\le m_0,
\end{array}
\right.
\end{equation}
then the Cauchy problem \eqref{ns}-\eqref{initial} admits
a unique global solution $(v,u,\t)$ satisfying
$$
(v-V,u-U,\theta-\Theta)(x,t)\in C\big((0,+\infty);H^1(\mathbb{R})\big);
$$
$$
(v-V)_x(x,t)\in L^2(0,+\infty;L^2\big(\mathbb{R})\big);
$$
$$
(u-U,\t-\Theta)_x(x,t)\in L^2\big(0,+\infty;H^1(\mathbb{R})\big),
$$
and
\begin{equation}\label{wt-2}
\lim_{t\rightarrow+\infty}\sup_{x\in\mathbb{R}}\left|(v-V,u-U,\t-\Theta)(x,t)\right|=0,
\end{equation}
where the $(v_-^r,u_-^r,\t_-^r)(x,t)$ and $(v_+^r,u_+^r,\t_+^r)(x,t)$ are the 1-rarefaction and
3-rarefaction waves uniquely determined by \eqref{euler}, \eqref{Riemann2}, respectively.
\end{theorem}

\begin{remark}
By $(iv)$ of Lemma \ref{rare-pro}, Theorem \ref{theorem2} implies
\begin{equation*}\label{behavior-2}
\lim_{t\rightarrow+\infty}\sup_{x\in\mathbb{R}}
\left(
\begin{array}{l}
|(v-v_-^r-V^{cd}-v_+^r+v_-^m+v_+^m)(x,t)|\\[2mm]
\quad\quad|(u-U^{cd}-u_-^r-u_+^r)(x,t)|\\[2mm]
|(\t-\t_-^r-\T^{cd}-\t_+^r+\t_-^m+\t_+^m)(x,t)|\\[2mm]
\end{array}
\right)=0.
\end{equation*}
\end{remark}

We now explain the main strategy of this paper. It is noted that in \cite{Hong,huang-zhao,nishi-yang-zhao},
the smallness of $\gamma-1$ is used to control the lower and upper bound of the absolute temperature $\t$.
To remove the smallness condition of $\gamma-1$, the key point is to derive the uniform bound of $\t$, which is
also closely related to the uniform bound of the specific volume $v$.
Motivated by \cite{jiang,kazhi-she,L-L}, we first obtain the basic energy estimate (see Lemma \ref{basic-lemma}),
which is independent of the time $t$,  with the help of the new estimates on the heat kernel developed in \cite{Huang-Li-matsumura},
provided the strengths of the waves are suitable small.
It should be emphasized that the basic energy estimate is nontrivially  obtained, while it is trivial for the case of small initial perturbation or the far-field condition being a constant one $(\bar{v},\bar{u},\bar{\t})$.
In fact, we essentially use the structure of wave patterns to control the terms involving the derivative of perturbation around the wave patterns.
Secondly, the specific volume $v$ is shown uniformly bounded from below and above with respect to space and time through
delicate analysis based on the basic energy estimate and a cut-off technique. Finally, we manipulate some weighted estimates on the perturbation
around the wave patterns to derive the uniform bound of $\t$. We remark that  the underlying structures of viscous contact wave and rarefaction waves are essentially used throughout the whole proof, and the idea may not be valid
for shock wave whose structure is quite different from those of
viscous contact wave and rarefaction waves.

\

This paper is organized as follows. In the next section, we collect some useful lemmas and
fundamental facts concerning the viscous contact wave as well as rarefaction waves.
The main proof of Theorem \ref{theorem1} and \ref{theorem2} are completed in Section 3 and 4, respectively.

\

\emph{Notations.} Throughout this paper, generic positive
constants are denoted by $c$ and $C$ without confusion.
For function spaces, $L^{p}(\Omega),1\leq p\leq
\infty$ denotes the usual Lebesgue space on $\Omega
\subset\mathbb{R}=(-\infty,\infty)$ with its norm given by
$$
\|f\|_{L^{p}(\Omega)}:=\left(\int_{\Omega}|f(x)|^{p}dx\right)^{\frac{1}{p}},
\quad 1\leq p<\infty, \quad  \parallel
f\parallel_{L^{\infty}(\Omega)}:=\mbox{ess.sup}_{\Omega} |f(x)|.
$$
$H^{k}(\Omega)$ denotes the $k^{th}$ order Sobolev space with its
norm
$$
\|f\|_{H^{k}(\Omega)}:=\left(\sum ^{k}_{j=0}
\parallel \partial^{j}_{x}f\parallel^{2}(\Omega)\right)^{\frac{1}{2}}, \quad \mathrm{when} \parallel
\cdot
\parallel=\parallel \cdot
\parallel_{L^{2}(\Omega)}.
$$
The domain $\Omega$ will be often abbreviated without confusion.


\

\section{Preliminaries}
\setcounter{equation}{0}

The properties of the viscous contact wave $(V,U,\Theta)$ defined by \eqref{contact} are useful  in the following sections.
\begin{lemma}\label{decay}
Assume that $\delta=|\t_+-\t_-|\leq \d_0$ for a small positive constant $\d_0$. Then the viscous contact wave $(V,U,\Theta)$ defined by \eqref{contact} has the following properties:
\begin{itemize}
\item [(1)]
$$
|V-v_{\pm}|+|\Theta-\t_{\pm}|\leq O(1)\d e^{-\frac{c_1x^2}{1+t}},
$$
\item [(2)]
$$
|\partial^k_x V|+|\partial^{k-1}_x U|+|\partial^k_x\Theta|\leq
O(1)\delta(1+t)^{-\frac{k}{2}}e^{-\frac{c_1x^2}{1+t}}, \quad k\geq 1.
$$
\end{itemize}

\end{lemma}
Therefore, we have
\begin{equation}
\di \widetilde R_1=O(1)\d(1+t)^{-\frac{3}{2}}e^{-\frac{c_1x^2}{1+t}},
\quad \widetilde R_2=O(1)\d(1+t)^{-2}e^{-\frac{c_1x^2}{1+t}}.
\end{equation}

The following two lemmas  play important roles to obtain the basic energy
estimate, the proofs can be found in \cite{Huang-Li-matsumura}, we omit them for brevity.
\begin{lemma}\label{hlm}
For $0<T\leq+\infty$, suppose that $h(x,t)$ satisfies
$$
h\in L^{\infty}(0,T;L^2(\mathbb{R})),\quad h_x\in L^{2}(0,T;L^2(\mathbb{R})),\quad h_t\in L^{2}(0,T;H^{-1}(\mathbb{R})).
$$
Then
\begin{equation}
\int_0^T\int h^2w^2dxdt\leq 4\pi\|h(0)\|^2+4\pi\alpha^{-1}\int_0^T\|h_x\|^2dt
+8\alpha\int_0^T<h_t,hg^2>_{H^{-1}\times H^1}dt
\end{equation}
for $\alpha>0$, and
$$
w(x,t)=(1+t)^{-\frac{1}{2}}\exp\left(-\frac{\alpha x^2}{1+t}\right),\quad
g(x,t)=\int_{-\infty}^x w(y,t)dy.
$$
\end{lemma}

\begin{lemma}\label{lemma5}
For $\alpha\in(0,\frac{c_1}{4}]$ and $w$ defined in Lemma \ref{hlm}, there exists some positive
constant $C$ depending on $\alpha$, such that the following estimate holds
\begin{equation}\label{important}
\int_0^t\int(\phi^2+\psi^2+\zeta^2)w^2 dxds\leq C\left(
1+\int_0^t\int (\phi_x^2+\psi_x^2+\zeta_x^2)dxds
\right).
\end{equation}

\end{lemma}

Next, we state the following properties of the solution to the problem \eqref{burgers}
due to \cite{MN-86}.

\begin{lemma}\label{bur-pro}
For given $w_l\in\mathbb{R}$ and $\bar w>0$, let $w_r\in\{0<\tilde w\triangleq w-w_l<\bar w\}$.
Then the problem \eqref{burgers} has a unique smooth global solution in time satisfying the following properties.
\begin{itemize}
\item[(i)] $w_l<w(x,t)<w_r$, $w_x>0$ $(x\in\mathbb{R},t>0)$.
\item[(ii)] For $p\in[1,\infty]$, there exists some positive constant $C=C(p,w_l,\bar w)$ such that
for $\tilde w\geqq0$ and $t\geqq0$,
$$
\|w_x(t)\|_{L^p}\leq C\min\{\tilde w,\tilde w^{1/p}t^{-1+1/p}\},\quad
\|w_{xx}(t)\|_{L^p}\leq C\min\{\tilde w,t^{-1}\}.
$$
\item[(iii)] If $w_l>0$, for any $(x,t)\in(-\infty,0]\times[0,\infty)$,
$$
|w(x,t)-w_l|\leq\tilde we^{-2(|x|+w_lt)},\quad
|w_x(x,t)|\leq2\tilde we^{-2(|x|+w_lt)}.
$$
\item[(iv)] If $w_r<0$, for any $(x,t)\in[0,\infty)\times[0,\infty)$,
$$
|w(x,t)-w_r|\leq\tilde we^{-2(x+|w_r|t)},\quad
|w_x(x,t)|\leq2\tilde we^{-2(x+|w_r|t)}.
$$
\item[(v)] For the Riemann solution $w^r(x/t)$ of the scalar equation \eqref{burgers} with the Riemann
initial data
\begin{equation*}
w(x,0)=\left\{
\begin{array}{ll}
w_l,&\quad x<0,\\
w_r,&\quad x>0,
\end{array}
\right.
\end{equation*}
we have
$$
\lim_{t\rightarrow+\infty}\sup_{x\in\mathbb{R}}|w(x,t)-w^r(x/t)|=0.
$$
\end{itemize}

\end{lemma}

Finally, we divide $\mathbb{R}\times(0,t)$ into
three parts, that is $\mathbb{R}\times(0,t)=\Omega_-\cup\Omega_c\cup\Omega_+$ with
$$
\Omega_{\pm}=\big\{(x,t)\big|\pm2x>\pm\lambda_{\pm}(v_{\pm}^m,s_{\pm})t\big\},
$$
and
$$
\Omega_{c}=\big\{(x,t)\big|\lambda_-(v_-^m,s_-)t\leq2x\leq\lambda_{+}(v_{+}^m,s_{+})t\big\}.
$$
Then Lemma \ref{decay} and Lemma \ref{bur-pro}  lead to

\begin{lemma}\label{rare-pro}
For any given left end state $(v_-,u_-,\t_-)$, we assume that \eqref{zheng}
and \eqref{same-order} hold. Then the smooth rarefaction waves $(V_{\pm}^r,U_{\pm}^r,\T_{\pm}^r)$
constructed in \eqref{appro-rare} and the viscous contact wave $(V^{cd},U^{cd},\T^{cd})$ constructed in \eqref{contact} satisfying
the following:
\begin{itemize}
\item[(i)] $(U_{\pm}^r)_x\geq0$,~$(x\in\mathbb{R},t>0)$.
\item[(ii)] For $p\in[1,\infty]$, there exists a positive constant $C=C(v_-,u_-,\t_-,\d)$
such that for $\d$ satisfying \eqref{same-order},
$$
\|\big((V_{\pm}^r)_x,(U_{\pm}^r)_x,(\T_{\pm}^r)_x\big)(t)\|_{L^p}
\leq C\min\Big\{\d,\d^{1/p}t^{-1+1/p}\Big\}
$$
and
$$
\|\big((V_{\pm}^r)_{xx},(U_{\pm}^r)_{xx},(\T_{\pm}^r)_{xx}\big)(t)\|_{L^p}
\leq C\min\Big\{\d,t^{-1}\Big\}.
$$
\item[(iii)] There exists a positive constant $C=C(v_-,u_-,\t_-,\d)$
such that for
$$
c_0=\frac{1}{10}\min\Big\{|\lambda_-(v_-^m,s_-)|,\lambda_+(v_+^m,s_+),c_1\lambda_-^2(v_-^m,s_-),
c_1\lambda_+^2(v_+^m,s_+),1\Big\},
$$
we have in $\Omega_c$
$$
(U_{\pm}^r)_x+|(V_{\pm}^r)_x|+|V_{\pm}^r-v_{\pm}^m|+|(\T_{\pm}^r)_x|+|\T_{\pm}^r-\t_{\pm}^m|
\leq C\d e^{-c_0(|x|+t)},
$$
and in $\Omega_{\mp}$
$$
\left\{
\begin{array}{ll}
|V^{cd}-v_{\mp}^m|+|V^{cd}_x|+|\T^{cd}-\t_{\mp}^m|+|U^{cd}_x|+|\T^{cd}_x|
\leq C\d e^{-c_0(|x|+t)},\\
(U_{\pm}^r)_x+|(V_{\pm}^r)_x|+|V_{\pm}^r-v_{\pm}^m|+|(\T_{\pm}^r)_x|+|\T_{\pm}^r-\t_{\pm}^m|
\leq C\d e^{-c_0(|x|+t)}.
\end{array}
\right.
$$
\item[(iv)] For the rarefaction waves $(v_{\pm}^r,u_{\pm}^r,\t_{\pm}^r)(x/t)$ determined by
\eqref{euler}\eqref{Riemann2}, it holds
$$
\lim_{t\rightarrow+\infty}\sup_{x\in\mathbb{R}}
\big|(V_{\pm}^r,U_{\pm}^r,\T_{\pm}^r,)(x,t)-(v_{\pm}^r,u_{\pm}^r,\t_{\pm}^r,)(x/t)\big|=0.
$$
\end{itemize}

\end{lemma}


\section{Proof of Theorem \ref{theorem1}}
\setcounter{equation}{0}

Substituting \eqref{contact equation} into \eqref{ns}, \eqref{initial} yields
\begin{equation}\label{perturb}
\left\{
\begin{array}{ll}
\di \phi_t-\psi_x=0,\\
\di \psi_t+\left(p-p_+\right)_x=\mu\left(\frac{u_x}{v}-\frac{U_x}{V}\right)_x-\widetilde R_1,\\
\di c_{\nu}\zeta_t+pu_x-p_+U_x=\kappa\left(\frac{\t_x}{v}-\frac{\Theta_x}{V}\right)_x
+\mu\left(\frac{u_x^2}{v}-\frac{U^2_x}{V}\right)-\widetilde R_2,\\
\di (\phi,\psi,\zeta)(x,0)=(\phi_0,\psi_0,\zeta_0)(x),\quad x\in \mathbb{R}.
\end{array}
\right.
\end{equation}
We shall prove Theorem \ref{theorem1} by the local existence and the a priori estimate.
We look for the solution $(\phi,\psi,\zeta)$ in the solution space $X([0,+\infty))$,
\begin{equation*}
\begin{array}{l}
\di X([0,T])=\Big\{(\phi,\psi,\zeta)\Big| v, ~\theta \geq M^{-1},
 \quad \sup_{0\leq t\leq T}\|(\phi,\psi,\zeta)\|_{H^1}\leq M
 \Big\}
\end{array}
\end{equation*}
for some $0<T\leq +\infty$, where the constants $M$ will be determined later.
Since the local existence of the solution is well known (for example, see \cite{huang-ma-shi}),
to prove the global existence part of Theorem \ref{theorem1}, we only need to establish the
following  a priori estimates.

\begin{proposition}\label{prop}
\textbf{(A priori estimates)} Assume that the conditions of Theorem \ref{theorem1} hold, then there exists
a positive constant $\d_0$ such that if $\d< \d_0$ and $(\phi,\psi,\zeta)\in X([0,T])$,
\begin{equation}
\begin{array}{l}\label{p-1}
\di \sup_{0\leq t\leq T}\|(\phi,\psi,\zeta)(t)\|_{H^1}^2+\int_0^T(\|\phi_x\|^2+\|(\psi_x,\zeta_x)\|^2_{H^1})ds
\leq C_0,
\end{array}
\end{equation}
where $C_0$ denotes a constant depending only on $\mu$, $\kappa$, $R$, $c_{\nu}$, $v_{\pm}$, $u_{\pm}$, $\t_{\pm}$ and  $m_0$.

\end{proposition}

Once Proposition \ref{prop} is proved, we can extend the unique local solution $(u,v,\t)$ which can be
obtained as in \cite{huang-ma-shi} to $T=\infty$. Estimate \eqref{p-1}
and the equations \eqref{perturb} (respectively \eqref{perturb2}) imply that
\begin{equation}
\di\int_0^{+\infty}\left(\|(\phi_x, \psi_x,\zeta_x)(t)\|^2+\left|\frac{d}{dt}\|(\phi_x, \psi_x,\zeta_x)(t)\|^2\right|\right)dt< \infty,
\end{equation}
which, together with \eqref{p-1} and the Sobolev's inequality, easily leads to the large time behavior of the solutions,
that is, \eqref{behavior-1} (resp. \eqref{wt-2}).

Proposition \ref{prop} will be finished by the following lemmas. First, we give the basic energy estimate,
which is nontrivially obtained, compared with the case of small initial perturbation or the far-field condition
being a constant one $(\bar v,\bar u,\bar\t)$.

\begin{lemma}\label{basic-lemma}
There exist some positive constant $C_0$  and $\d_0$ such that if $\d<\d_0$, it holds that
\begin{equation}\label{basic}
\begin{array}{ll}
\di \Big\|\left(\psi,\sqrt{\Phi\left(\frac{v}{V}\right)},\sqrt{\Phi\left(\frac{\t}{\Theta}\right)}\right)(t)\Big\|^2
+\int_0^t\int\left(\frac{\psi_x^2}{\t v}+\frac{\zeta_x^2}{\t^2 v}\right)dxds
 \leq C_0.
\end{array}
\end{equation}

\end{lemma}

\textbf{Proof}: The proof of the Lemma \ref{basic-lemma} consists of two steps.

\underline{ Step 1. }\quad Similar to \cite{huang-ma-shi}, multiplying $\eqref{ns}_1$ by
$-R\Theta(v^{-1}-V^{-1})$, $\eqref{ns}_2$ by $\psi$ and $\eqref{ns}_3$ by $\zeta\theta^{-1}$, then
adding the resulting equations together, we get
\begin{equation}
\begin{array}{ll}
\di\left(\frac{\psi^2}{2}+R\Theta\Phi\left(\frac{v}{V}\right)+c_{\nu}\Theta\Phi\left(\frac{\theta}{\Theta}\right)\right)_t
+\frac{\mu\Theta}{\theta v}\psi_x^2+\frac{\kappa\Theta}{\theta^2 v}\zeta_x^2
+H_{x}+Q=-\psi \widetilde R_1-\frac{\zeta}{\theta}\widetilde R_2
\end{array}
\label{ba1}
\end{equation}
with
$$
\Phi(z)=z-\ln z-1, \quad z>0
$$
and
\begin{equation}\label{H}
\di H=(p-p_+)\psi-\mu\left(\frac{u_x}{v}-\frac{U_x}{V}\right)\psi-\frac{\k\zeta}{\theta}\left(\frac{\t_x}{v}-\frac{\Theta_x}{V}\right),
\end{equation}
\begin{equation}
\begin{array}{ll}
\di Q&\di=p_+\Phi\left(\frac{V}{v}\right)U_x+\frac{p_+}{\gamma-1}\Phi\left(\frac{\Theta}{\theta}\right)U_x
+\mu\left(\frac{1}{v}-\frac{1}{V}\right)\psi_xU_x\\[3mm]
&\di -\frac{\zeta}{\theta}(p_+-p)U_x-\frac{\k\Theta_x}{\theta^2v}\zeta\zeta_x
-\frac{\k\Theta\Theta_x}{\theta^2vV}\zeta_x\phi+\frac{\k\Theta_x^2}{\t^2vV}\zeta\phi\\[3mm]
&\di-\frac{2\mu U_x}{\theta v}\psi_x\zeta+\frac{\mu U_x^2}{\theta vV}\zeta\phi.
\end{array}
\end{equation}
Since
\begin{equation}
\di Q\leq \frac{\mu\Theta}{4\theta v}\psi_x^2+\frac{\k\Theta}{4\theta^2v}\zeta_x^2
+C(M)(\phi^2+\zeta^2)(|U_x|+\Theta_x^2),
\end{equation}
where $C(M)$ denotes a constant depending on $M$. Recalling Lemma \ref{decay}, we have
\begin{equation}
\begin{array}{ll}
\di \left|\int_0^t\int\widetilde R_1\psi dxds\right|&\di\leq O(1)\d \int_0^t\int (1+s)^{-3/2}e^{-\frac{c_1x^2}{1+s}}|\psi|dxds\\[3mm]
&\di \leq O(1)\d\int_0^t(1+s)^{-5/4}\|\psi\|ds\\[3mm]
&\di \leq O(1)\d\int_0^t(1+s)^{-5/4}\|\psi\|^2ds+O(1)\d
\end{array}
\end{equation}
and
\begin{equation}
\begin{array}{ll}
\di \left|\int_0^t\int\widetilde R_2\frac{\zeta}{\t} dxds\right|
\leq C(M)\d\int_0^t(1+s)^{-7/4}\left\|\sqrt{\Phi\left(\frac{\t}{\T}\right)}\right\|^2ds+C(M)\d.
\end{array}
\end{equation}
Then integrating \eqref{ba1} over $\mathbb{R}\times(0,t)$, choosing $\a=\frac{c_1}{4}$ in Lemma \ref{hlm}
and $\d$ suitable small, it follows from Lemma \ref{hlm}-\ref{lemma5} and Gronwall's inequality that
\begin{equation}\label{ba2}
\begin{array}{ll}
\di \Big\|\left(\psi,\sqrt{\Phi\left(\frac{v}{V}\right)},\sqrt{\Phi\left(\frac{\t}{\Theta}\right)}\right)(t)\Big\|^2
+\int_0^t\int\left(\frac{\psi_x^2}{\t v}+\frac{\zeta_x^2}{\t^2 v}\right)dxds\\[3mm]
\di
\leq C_0+ C(M)\d \int_0^t\int(1+s)^{-1}(\phi^2+\zeta^2)e^{-\frac{c_1x^2}{1+s}}dxds\\
\di \leq C_0+ C(M)\d \int_0^t\int\frac{\t\phi_x^2}{v^3}dxds.
\end{array}
\end{equation}

\underline{ Step 2. }\quad  Following \cite{MN-92}, we introduce a new
variable $\di\tilde v=\frac{v}{V}$. Then $\eqref{perturb}_2$ can be rewritten by the
new variable as
\begin{equation}\label{newperturb}
\di \left(\mu\frac{\tilde v_x}{\tilde v}-\psi\right)_t-p_x=\widetilde R_1.
\end{equation}
Multiplying \eqref{newperturb} by $\di\frac{\tilde v_x}{\tilde v}$, we have
\begin{equation}\label{phi-x1}
\begin{array}{ll}
\di \left(\frac{\mu}{2}\left(\frac{\tilde v_x}{\tilde v}\right)^2-\psi\frac{\tilde v_x}{\tilde v}\right)_t
+\left(\psi\frac{\tilde v_t}{\tilde v}\right)_x+\frac{R\t}{v}\left(\frac{\tilde v_x}{\tilde v}\right)^2
-\frac{R}{v}\zeta_x\frac{\tilde v_x}{\tilde v}\\[3mm]
\di\quad +\frac{R\t}{v}\left(\frac{1}{\Theta}-\frac{1}{\t}\right)\Theta_x\frac{\tilde v_x}{\tilde v}
=\frac{\psi_x^2}{v}+\psi_xU_x\left(\frac{1}{v}-\frac{1}{V}\right)+\widetilde R_1\frac{\tilde v_x}{\tilde v}.
\end{array}
\end{equation}
The Cauchy's inequality yields that
\begin{equation}
\begin{array}{ll}
\di \left|\frac{R}{v}\zeta_x\frac{\tilde v_x}{\tilde v}\right|
+\left|\psi_xU_x\left(\frac{1}{v}-\frac{1}{V}\right)\right|+\frac{\psi_x^2}{v}\\[3mm]
\di \leq \frac{R\t}{4v}\left(\frac{\tilde v_x}{\tilde v}\right)^2+
C(M)\left(\frac{\zeta_x^2}{\t^2 v}+\frac{\psi_x^2}{\t v}\right)
+C(M)\phi^2U_x^2,
\end{array}
\end{equation}
and
\begin{equation}
\begin{array}{ll}
\di \left|\frac{R\t}{v}\left(\frac{1}{\Theta}-\frac{1}{\t}\right)\Theta_x\frac{\tilde v_x}{\tilde v}\right|
+\left|\widetilde R_1\frac{\tilde v_x}{\tilde v}\right|\\[3mm]
\di\leq \frac{R\t}{4v}\left(\frac{\tilde v_x}{\tilde v}\right)^2
+C(M)(\zeta^2\Theta_x^2+\widetilde R_1^2).
\end{array}
\end{equation}
Note that
\begin{equation*}
\begin{array}{ll}
\di \frac{\phi_x^2}{2v^2}-C(M)\phi^2\Theta_x^2\leq \left(\frac{\tilde v_x}{\tilde v}\right)^2
\leq \frac{\phi_x^2}{v^2}+C(M)\phi^2\Theta_x^2.
\end{array}
\end{equation*}
Integrating \eqref{phi-x1} over $\mathbb{R}\times(0,t)$, we have
\begin{equation}
\begin{array}{ll}
\di\int\frac{\phi_x^2}{v^2}dx+\int_0^t\int\frac{\t\phi_x^2}{v^3} dxds\leq C_0+C\|\psi\|^2\\[3mm]
\di\quad+C(M)\d^2\left\|\sqrt{\Phi\left(\frac{v}{V}\right)}\right\|^2
+C(M)\int_0^t\int\left(\frac{\psi_x^2}{\t v}+\frac{\zeta_x^2}{\t^2 v}\right) dxds\\[3mm]
\di\quad +C(M)\int_0^t\int(\phi^2+\zeta^2)(U_x^2+\Theta_x^2) dxds.
\end{array}
\end{equation}
By Lemma \ref{lemma5}, \eqref{ba2} and choosing $\d$ suitable small, we have
\begin{equation}\label{phi-x3}
\begin{array}{ll}
\di \int\frac{\phi_x^2}{v^2}dx+\int_0^t\int\frac{\theta\phi_x^2}{v^3}dxds
 \leq C_0+C(M)\int_0^t\int\left(\frac{\psi_x^2}{\t v}+\frac{\zeta_x^2}{\theta^2v}\right) dxds.
\end{array}
\end{equation}
Then the proof of Lemma \ref{basic-lemma} is completed by substituting \eqref{phi-x3} into \eqref{ba2},
and choosing $\d$ suitable small.
\hfill $\Box$

In Lemma \ref{basic-lemma}, the smallness of $\delta$ is used to guarantee that the basic energy (\ref{basic}) is only bounded by the initial data. Based on
the basic energy estimate, we shall show the specific volume $v$ and the absolute temperature $\t$ are uniformly bounded from below and above, which in turn decides
how small for $\delta$. That is why need the initial condition (\ref{data}). To prove Theorem \ref{theorem1}, we first try to get the uniform bound of $v(x,t)$. We have
\begin{lemma}
Let $\a_1$, $\a_2$ be the two positive roots of the equation $y-\ln y-1=C_0$ and the constant $C_0$ be the same in \eqref{basic}. Then
\begin{equation}\label{root}
\di \a_1\leq \int_k^{k+1}\tilde v(x,t)dx, \quad\int_k^{k+1}\tilde \t(x,t)dx\leq \a_2,\quad t\geq 0,
\end{equation}
and for each $t\geq 0$ there are points $a_k(t)$, $b_k(t)\in[k,k+1]$ such that
\begin{equation}\label{root2}
\alpha_1\leq \tilde v(a_k(t),t), \tilde \t(b_k(t),t)\leq \alpha_2,\quad t\geq 0,
\end{equation}
where $k=0,\pm1,\pm2, \cdots$.
\end{lemma}

\textbf{Proof}: From \eqref{basic}, we see that
\begin{equation}\label{4.1}
\di \int_k^{k+1}(\tilde v(x,t)-\ln\tilde v(x,t)-1)dx, \quad\int_k^{k+1}(\tilde \t(x,t)-\ln\tilde \t(x,t)-1)dx\leq C_0,
\end{equation}
where $\tilde v=\frac{v}{V}$, $\tilde \t=\frac{\t}{\Theta}$, and $k=0,\pm1,\pm2,\cdots$. If we apply Jessen's inequality to the convex function
$y-\ln y-1=C_0$, we obtain
$$
\di \int_k^{k+1}\tilde v(x,t)dx-\ln\int_k^{k+1}\tilde v(x,t)dx-1,
\quad\int_k^{k+1}\tilde \t(x,t)dx-\ln\int_k^{k+1}\tilde \t(x,t)dx-1\leq C_0,
$$
which gives
$$
\di \a_1\leq \int_k^{k+1}\tilde v(x,t)dx, \quad\int_k^{k+1}\tilde \t(x,t)dx\leq \a_2,
$$
where $\alpha_1$, $\alpha_2$ are two positive roots of the equation $y-\ln y-1=C_0$.
Moreover, in view of mean value theorem, for each $t\geq 0$, there are points
$a_k(t)$, $b_k(t)\in [k,k+1]$ such that
\begin{equation}
\di 0<\alpha_1\leq \tilde v(a_k(t),t), \tilde \t(b_k(t),t)\leq \alpha_2,\quad t\geq 0.
\end{equation}

\hfill $\Box$

The following Lemma can be found in Jiang \cite{jiang} Lemma 2.3.
\begin{lemma}\label{jiang-lemma}
For each $x\in[k.k+1]$, $k=0,\pm1,\pm2,\cdots$, it follows from \eqref{ns}$_2$ that
\begin{equation}\label{v-repre}
\di v(x,t)=B(x,t)Y(t)+\frac{R}{\mu}\int_0^t\frac{B(x,t)Y(t)}{B(x,s)Y(s)}\t(x,s) ds,
\end{equation}
where
\begin{equation}\label{B}
\di B(x,t)=v_0(x)\exp\left(\frac{1}{\mu}\int_x^{+\infty}\big(u_0(y)-u(y,t)\big)\beta(y)dy\right),
\end{equation}
\begin{equation}\label{Y}
\di Y(t)=\exp\left(\frac{1}{\mu}\int_0^t\int_{k+1}^{k+2}\s(y,s)dyds\right),
\end{equation}
\begin{equation}\label{sigma}
\di \sigma(x,t)=\left(\mu\frac{u_x}{v}-R\frac{\t}{v}\right)(x,t),
\end{equation}
and
\begin{equation}\label{beta}
\di \beta(x)=\left\{
\begin{array}{ll}
1,& x\leq k+1,\\
k+2-x, & k+1\leq x\leq k+2,\\
0, & x\geq k+2.
\end{array}
\right.
\end{equation}

\end{lemma}
By Cauchy's inequality and \eqref{basic}, we have
\begin{equation}\label{b-bound}
\di \underline{B}(C_0)\leq B(x,t)\leq \overline{B}(C_0), \quad \forall x\in[k,k+1],\quad t\geq 0,
\end{equation}
where $\underline{B}(C_0)$, $\overline{B}(C_0)$ are two constants depending on $C_0$.

\begin{lemma}\label{v-bound}
There are two positive constants $\underline{v}(C_0)$, $\bar{v}(C_0)$  such that
\begin{equation}\label{v}
\di \underline{v}(C_0) \leq v(x,t)\leq\bar{v}(C_0),\quad \forall x\in\mathbb{R},\quad t\geq0,
\end{equation}
where $\underline{v}(C_0)$, $\bar{v}(C_0)$  depending on $C_0$, independent of $x$, $t$.

\end{lemma}
\textbf{Proof}: From now on, we always assume that $\t_-<\t_+$ for convenient. So from the properties of
viscous contact wave, we have $\t_-<\T(x,t)<\t_+$ and $v_-<V(x,t)<v_+$.
For each $t\geq0$, there exists at least one point $x=x_{k+1}(t)\in[k+1,k+2]$ such that
$$
\di \inf_{x\in[k+1,k+2]}\tilde\t(x,t)=\tilde\t(x_{k+1}(t),t).
$$
By Cauchy's inequality, \eqref{basic}, \eqref{root}, and choosing $\d$ suitable small, we see that
\begin{equation}\label{4.2}
\begin{array}{ll}
\di\left|\int_s^t\int_{b_{k+1}(\tau)}^{x_{k+1}(\tau)}\frac{\tilde\t_y}{\tilde\t}(y,\tau) dyd\tau \right|=
\left|\int_s^t\int_{b_{k+1}(\tau)}^{x_{k+1}(\tau)} \left(\frac{\zeta_y}{\t}-\frac{\zeta\Theta_y}{\t\Theta}\right)dyd\tau \right|\\[3mm]
\di\quad\leq\int_s^t\left(\int_{k+1}^{k+2}\frac{\zeta_x^2}{\t^2v}dx\right)^{\frac{1}{2}}\left(\int_{k+1}^{k+2}vdx\right)^{\frac{1}{2}}d\tau
+\int_s^t\int_{k+1}^{k+2}\left|\frac{\zeta\Theta_x}{\t\Theta}\right|dxd\tau\\[3mm]
\di\quad\leq\int_s^t\left(\int_{k+1}^{k+2}\frac{\zeta_x^2}{\t^2v}dx\right)^{\frac{1}{2}}\big(v_+\a_2\big)^{\frac12}d\tau
+C(M)\int_s^t\left(\int_{k+1}^{k+2}\zeta^2\Theta_x^2 dx\right)^{\frac12}d\tau\\[3mm]
\di\quad \leq C\left(\int_s^t\int\frac{\zeta_x^2}{\t^2v} dyd\tau\right)^{\frac12}\sqrt{t-s}
+C(M)\left(\int_s^t\int\zeta^2\Theta_x^2 dxd\tau\right)^{\frac12}\sqrt{t-s}\\[3mm]
\di\quad \leq C_0\sqrt{t-s}.
\end{array}
\end{equation}
We apply Jessen's inequality to the convex function $e^x$, and utilize \eqref{basic}, \eqref{root2},
\eqref{4.2} to obtain that for $t\geq s\geq 0$,
\begin{equation}
\begin{array}{ll}
\di\int_s^t\inf_{x\in[k+1,k+2]}\tilde \t(\cdot,\tau) d\tau=
\int_s^t\tilde \t(x_{k+1}(\tau),\tau) d\tau=
\int_s^t\exp\left(\log\tilde\t(x_{k+1}(\tau),\tau)\right) d\tau\\[4mm]
\di\quad\geq(t-s)\exp\left(\frac{1}{t-s}\int_s^t\log\tilde\t(x_{k+1}(\tau),\tau) d\tau\right)\\[4mm]
\di\quad=(t-s)\exp\left(\frac{1}{t-s}\int_s^t\left[\log\frac{\tilde\t(x_{k+1}(\tau),\tau)}{\tilde\t(b_{k+1}(\tau),\tau)}
+\log\tilde\t(b_{k+1}(\tau),\tau)\right] d\tau\right)\\[4mm]
\di\quad=(t-s)\exp\left(\frac{1}{t-s}\int_s^t\left[\int_{b_{k+1}(\tau)}^{x_{k+1}(\tau)}\frac{\tilde\t_y}{\tilde\t} dy
+\log\tilde\t(b_{k+1}(\tau),\tau)\right]d\tau\right)\\[4mm]
\di\quad\geq(t-s)\exp\left(\log\a_1-\frac{1}{t-s}\left|\int_s^t\int_{b_{k+1}(\tau)}^{x_{k+1}(\tau)}\frac{\tilde\t_y}{\tilde\t}dyd\tau\right|\right)\\[4mm]
\di\quad\geq C(t-s)e^{-\frac{C_0}{\sqrt{t-s}}}.
\end{array}
\end{equation}
Noticing that $\t_-\leq\Theta\leq\t_+$, so we have
\begin{equation}\label{int-theta}
-\int_s^t\inf_{x\in[k+1,k+2]}\t(\cdot,\tau) d\tau\leq
\left\{
\begin{array}{ll}
0,&0\leq t-s\leq 1,\\
-C_0(t-s),&t-s\geq 1.
\end{array}
\right.
\end{equation}
Applying Cauchy's inequality and Jessen's inequality for the function $\frac{1}{x}(x>0)$, using
\eqref{basic}, \eqref{int-theta}, and noting that
\begin{equation}
\di\left\{
\begin{array}{ll}
C_0, &0\leq t-s\leq 1,\\
C_0-\frac{t-s}{C_0}, & t-s\geq 1
\end{array}
\right.
\leq C_0-\frac{t-s}{C_0},\quad t\geq s\geq 0,
\end{equation}
we obtain
\begin{equation}\label{int-sigma}
\begin{array}{ll}
\di\int_s^t\int_{k+1}^{k+2}\sigma(x,\tau)dxd\tau=
\int_s^t\int_{k+1}^{k+2}\left(\mu\frac{\psi_x}{v}-R\frac{\t}{v}+\mu\frac{U_x}{v}\right)(x,\tau)dxd\tau\\[3mm]
\di\quad\leq C\int_s^t\int_{k+1}^{k+2}\frac{\psi_x^2}{\t v}dxd\tau-\frac{R}{2}\int_s^t\int_{k+1}^{k+2}\frac{\t}{v}dxd\tau
+\mu\int_s^t\int_{k+1}^{k+2}\frac{U_x}{v}dxd\tau\\[3mm]
\di\quad\leq C_0-\frac{R}{2}\int_s^t\inf_{x\in[k+1,k+2]}\t\left(\int_{k+1}^{k+2}v^{-1}dx\right)d\tau
+C(M)\left(\int_s^t\int_{k+1}^{k+2} U_x^2 dxd\tau\right)^{\frac{1}{2}}(t-s)^{\frac12}\\[3mm]
\di\quad\leq C_0-\frac{R}{2}\int_s^t\inf_{x\in[k+1,k+2]}\t\left(\int_{k+1}^{k+2}vdx\right)^{-1}d\tau
+C(M)\d(t-s)^{\frac12}\\[3mm]
\di\quad\leq C_0-\frac{R}{2\a_2v_+}\int_s^t\inf_{x\in[k+1,k+2]}\t d\tau
+C(t-s)^{\frac12}\\[3mm]
\di\quad\leq C_0-\frac{t-s}{C_0}.
\end{array}
\end{equation}
It follows from  the definition of $Y(t)$ and \eqref{int-sigma} that
\begin{equation}\label{Y-e}
\di 0\leq Y(t)\leq C_0e^{-t/C_0},\quad
\frac{Y(t)}{Y(s)}\leq C_0e^{-(t-s)/C_0},
\end{equation}
which, together with \eqref{v-repre} and \eqref{b-bound}, gives,
\begin{equation}\label{v-repre2}
\di v(x,t)\leq C_0+C_0\int_0^t\t(x,s)e^{-(t-s)/C_0} ds.
\end{equation}
On the other hand, we have
\begin{equation}
\begin{array}{ll}
\di |\tilde \t^{\frac{1}{2}}(x,t)-\tilde\t^{\frac{1}{2}}(b_k(t),t)|
\leq\int_{k}^{k+1}\tilde\t^{-\frac{1}{2}}|\tilde\t_x|dx
\leq\int_k^{k+1}\left(\frac{\Theta}{\theta}\right)^{\frac{1}{2}}
\left(\left|\frac{\zeta_x}{\Theta}\right|+\left|\frac{\zeta\Theta_x}{\Theta^2}\right|\right)dx\\[3mm]
\di\quad\leq\t_{-}^{-\frac{1}{2}}\int_k^{k+1}\frac{|\zeta_x|}{\sqrt{\t}}dx
+C(M)\int_k^{k+1}|\zeta\Theta_x|dx\\[3mm]
\di\quad\leq\t_{-}^{-\frac{1}{2}}\left(\int_k^{k+1}\frac{\zeta_x^2}{\t^2v}dx\right)^{\frac{1}{2}}
\left(\int_k^{k+1}\t vdx\right)^{\frac{1}{2}}+C(M)\left(\int_k^{k+1}\zeta^2\Theta_x^2dx\right)^{\frac{1}{2}}\\[3mm]
\di\quad\leq\t_{-}^{-\frac{1}{2}}\left(\int_k^{k+1}\frac{\zeta_x^2}{\t^2v}dx\right)^{\frac{1}{2}}
\left(\int_k^{k+1}\t dx\right)^{\frac{1}{2}}\max_{x\in[k,k+1]}v(\cdot,t)^{\frac{1}{2}}
+C(M)\left(\int_k^{k+1}\zeta^2\Theta_x^2dx\right)^{\frac{1}{2}}\\[3mm]
\di\quad\leq\sqrt{\frac{\a_2\t_+}{\t_-}}\left(\int\frac{\zeta_x^2}{\t^2v}dx\right)^{\frac{1}{2}}\max_{x\in[k,k+1]}v(\cdot,t)^{\frac{1}{2}}
+C(M)\left(\int\zeta^2\Theta_x^2dx\right)^{\frac{1}{2}}\quad \mathrm{for}\quad x\in[k,k+1]
\end{array}
\end{equation}
and $k=0$, $\pm1$, $\pm2, \cdots$, which, along with \eqref{root2}, leads to
\begin{equation}\label{v-theta}
\begin{array}{ll}
\di \frac{\a_1\t_-}{3}-\a_2\frac{\t_+^2}{\t_-}\left(\int_{\mathbb{R}}\frac{\zeta_x^2}{\t^2v}\right)\max_{x\in\mathbb{R}}v(\cdot,t)
-C(M)\int_{\mathbb{R}}\zeta^2\Theta^2_x\\[3mm]
\di\quad\leq \theta(x,t)\leq 3\a_2\t_+ +3\a_2\frac{\t_+^2}{\t_-}\left(\int_{\mathbb{R}}\frac{\zeta_x^2}{\t^2v}\right)\max_{x\in\mathbb{R}}v(\cdot,t)
+C(M)\int_{\mathbb{R}}\zeta^2\Theta^2_x, \quad \forall x\in\mathbb{R}.
\end{array}
\end{equation}
Hence, substituting \eqref{v-theta} into \eqref{v-repre2}, applying Gronwall's inequality and \eqref{basic}, one has,
\begin{equation}\label{v-upper}
\di v(x,t)\leq C_0,\qquad \forall x\in\mathbb{R},\quad t\geq 0.
\end{equation}
Integrating \eqref{v-repre} over $[k,k+1]$ with respect to $x$, we obtain
\begin{equation}
\begin{array}{ll}
\di v_-\a_1&\di\leq C_0e^{-t/C_0}+C_0\int_0^t\frac{Y(t)}{Y(s)}\int_k^{k+1}\theta(x,s)dxds\\
&\di \leq C_0e^{-t/C_0}+C_0\int_0^t\frac{Y(t)}{Y(s)}ds.
\end{array}
\end{equation}
This directly yields that
\begin{equation}\label{Yt-Ys}
\di \int_0^t\frac{Y(t)}{Y(s)} ds\geq C_0-C_0e^{-t/C_0}.
\end{equation}
From \eqref{basic},\eqref{v-repre}, \eqref{v-theta}, \eqref{v-upper} and \eqref{Yt-Ys}, and choosing $\d$
suitable small, we have
\begin{equation}\label{wt-1}
\begin{array}{ll}
\di v(x,t)&\di\geq C_0\int_0^t\frac{Y(t)}{Y(s)}\theta(x,s)ds\\
&\di\geq C_0\int_0^t\frac{Y(t)}{Y(s)}ds-C_0\left(\int_0^{t/2}+\int_{t/2}^t\right)\frac{Y(t)}{Y(s)}\int\frac{\zeta_x^2}{\t^2v} dxds\\
&\di\qquad-C(C_0,M)\int_0^t\int\zeta^2\Theta_x^2 dxds\\
&\di\geq C_0-C_1(C_0)e^{-t/C_0}-C_0e^{-t/2C_0}\int_0^t\int\frac{\zeta_x^2}{\theta^2v} dxds
 -C_0\int_{t/2}^t\int\frac{\zeta_x^2}{\theta^2v}dxds\\
&\di\geq C_0/2, \qquad \forall x\in\mathbb{R}, ~t\geq T_0,
\end{array}
\end{equation}
where $C_1(C_0)$ is some positive constant depending on $C_0$, and $T_0$, $C_0$ are positive constants independent of $t$.

Next we consider the lower bound of $v(x,t)$ on $[0,T]$ for a positive constant $T>0$. From \cite{kazhikhov},
for any $x\in[k,k+1]$, $k=0,\pm1,\pm2,\cdots$, it holds that
\begin{equation}\label{v-ka}
\di v(x,t)=\frac{1}{\widetilde Y(t)\widetilde B(x,t)}
\left(v_0(x)+\frac{R}{\mu}\int_0^t\widetilde Y(s)\widetilde B(x,s)\t(x,s) ds\right),
\end{equation}
where
$$
\di\widetilde Y(t)=\frac{v_0(a_k(t))}{v(a_k(t),t)}\exp\left(\frac{R}{\mu}\int_0^t\frac{\t}{v}(a_k(t),s)ds\right)
$$
and
$$
\di\widetilde B(x,t)=\exp\left(\frac{1}{\mu}\int_{a_k(t)}^x\Big(u_0(y)-u(y,t)\Big)dy\right)
$$
with $a_{k}(t)$ is the same as in \eqref{root2}.
It follows from \eqref{v-ka} that
\begin{equation}\label{v-ka2}
\di \widetilde Y(t)v(x,t)=\frac{1}{\widetilde B(x,t)}
\left(v_0(x)+\frac{R}{\mu}\int_0^t\widetilde Y(s)\widetilde B(x,s)\t(x,s) ds\right).
\end{equation}
Integrating \eqref{v-ka2} over $[k,k+1]$ with respect to $x$, we obtain
\begin{equation}
\begin{array}{ll}
\di \a_1v_- \widetilde Y(t)&\di\leq C_0+C_0\int_0^t\widetilde  Y(s)\int_k^{k+1}\t(x,s) dxds\\
&\di\leq C_0+C_0\int_0^t\widetilde  Y(s)ds,
\end{array}
\end{equation}
which, together with the Gronwall's inequality, yields
$$
\widetilde  Y(t)\leq C(C_0,T).
$$
Then from \eqref{v-ka}, one has
\begin{equation*}
\di v(x,t)\geq \frac{v_0(x)}{\widetilde Y(t)\widetilde B(x,t)}\geq C(C_0,T),
\qquad \forall x\in \mathbb{R},
\end{equation*}
which, together with \eqref{wt-1} and \eqref{v-upper}, completes the proof of Lemma \ref{v-bound}.
\hfill $\Box$

\

Motivated by \cite{L-L}, we shall show the uniform bound of the absolute temperature $\t$
from below and above with respect to space and time. We have

\begin{lemma}\label{lemma-ll-1}
There exists some positive constants $C_0$ such that for any given $T>0$,
\begin{equation}\label{new-nergy}
\di\sup_{0\leq t\leq T}\int(\zeta^2+\psi^4)dx+\int_0^T\int((\t+\psi^2)\psi_x^2+\zeta_x^2) dxdt
\leq C_0.
\end{equation}
\end{lemma}
\textbf{Proof}: The proof of the Lemma \ref{lemma-ll-1} consists of the following steps.

\underline{ Step 1. }\quad First, for $t\geq0$, and $a>1$, denoting
$$
\Omega_a(t)\triangleq \left\{x\in\mathbb{R}\Bigg|\frac{\t}{\Theta}(x,t)>a\right\}
=\{x\in\mathbb{R}|\zeta(x,t)>(a-1)\Theta(x,t)\}.
$$
We derive from \eqref{basic} that $\Omega_a$ is bounded since
\begin{equation}\label{omega-bound}
\di a|\Omega_a|<\sup_{0\leq t\leq T}\int_{\Omega_a}\frac{\t}{\Theta}dx
\leq C(a)\sup_{0\leq t\leq T}\int_{\mathbb{R}}\Phi\left(\frac{\t}{\Theta}\right)dx\leq C(a,C_0).
\end{equation}

Next, multiplying $\eqref{perturb}_3$ by $(\zeta-\Theta)_+=\max\{\zeta-\Theta,0\}$, then integrating the resulted equation
over $\mathbb{R}\times[0,t]$, one has
\begin{equation}\label{5.1}
\begin{array}{ll}
\di\frac{c_{\nu}}{2}\int(\zeta-\Theta)_+^2dx+\k\int_0^t\int_{\Omega_2}\frac{\zeta_x^2}{v}dxds
=\frac{c_{\nu}}{2}\int(\zeta_0(x)-\Theta(x,0))_+^2dx\\[2mm]
\di\quad-\int_0^t\int\frac{R\zeta+R\Theta}{v}\psi_x(\zeta-\Theta)_+ dxds
-\int_0^t\int\frac{R\zeta-p_+\phi}{v}U_x(\zeta-\Theta)_+ dxds\\[2mm]
\di\quad+\k\int_0^t\int_{\Omega_2}\frac{\zeta_x\Theta_x}{V}dxds
-\k\int_0^t\int_{\Omega_2}\frac{\phi\Theta_x^2}{vV}dxds
+\mu\int_0^t\int\frac{\psi_x^2}{v}(\zeta-\Theta)_+ dxds\\[2mm]
\di\quad+2\mu\int_0^t\int\frac{\psi_xU_x}{v}(\zeta-\Theta)_+ dxds-\mu\int_0^t\int\frac{\phi U_x^2}{vV}(\zeta-\Theta)_+ dxds\\[2mm]
\di\quad-\int_0^t\int\widetilde R_2(\zeta-\Theta)_+ dxds-c_{\nu}\int_0^t\int\partial_t\Theta(\zeta-\Theta)_+ dxds.
\end{array}
\end{equation}
We multiply $\eqref{perturb}_2$ by $2\psi(\zeta-\Theta)_+$, and integrate the resulting equation over
$\mathbb{R}\times[0,t]$ to get
\begin{equation}\label{5.2}
\begin{array}{ll}
\di \int\psi^2(\zeta-\Theta)_+dx+2\mu\int_0^t\int\frac{\psi_x^2}{v}(\zeta-\Theta)_+ dxds
=\int\psi_0^2(x)(\zeta_0(x)-\Theta(x,0))_+dx\\[2mm]
\di\quad+2\int_0^t\int\frac{R\zeta-p_+\phi}{v}\psi_x(\zeta-\Theta)_+dxds
+2\int_0^t\int_{\Omega_2}\frac{R\zeta-p_+\phi}{v}\psi\zeta_xdxds\\[2mm]
\di\quad-2\int_0^t\int_{\Omega_2}\frac{R\zeta-p_+\phi}{v}\psi\Theta_xdxds
+2\mu\int_0^t\int\frac{\phi U_x}{vV}\psi_x(\zeta-\Theta)_+ dxds\\[2mm]
\di\quad-2\mu\int_0^t\int_{\Omega_2}\frac{\psi\psi_x\zeta_x}{v}dxds
+2\mu\int_0^t\int_{\Omega_2}\frac{\phi\psi U_x}{vV}\zeta_xdxds\\[2mm]
\di\quad+2\mu\int_0^t\int_{\Omega_2}\frac{\psi\psi_x\Theta_x}{v}dxds
-2\mu\int_0^t\int_{\Omega_2}\frac{\phi\psi}{vV}U_x\Theta_xdxds\\[2mm]
\di\quad-2\int_0^t\int\psi\widetilde R_1(\zeta-\Theta)_+ dxds
+\int_0^t\int_{\Omega_2}\psi^2\partial_t\zeta dxds-\int_0^t\int_{\Omega_2}\psi^2\partial_t\Theta dxds.
\end{array}
\end{equation}
Adding \eqref{5.2} into \eqref{5.1}, using $\eqref{perturb}_3$, we have
\begin{equation}
\begin{array}{ll}\label{long}
\di\int\left(\frac{c_{\nu}}{2}(\zeta-\Theta)_+^2+\psi^2(\zeta-\Theta)_+\right)dx
+\mu\int_0^t\int\frac{\psi_x^2}{v}(\zeta-\Theta)_+ dxds+\k\int_0^t\int_{\Omega_2}\frac{\zeta_x^2}{v}dxds\\[3mm]
\di=\int\left(\frac{c_{\nu}}{2}(\zeta_0(x)-\Theta(x,0))_+^2+\psi_0^2(x)(\zeta_0(x)-\Theta(x,0))_+\right)dx
+\int_0^t\int\frac{R\zeta-2p_+\phi-R\Theta}{v}\psi_x(\zeta-\Theta)_+ dxds\\[3mm]
\di-\int_0^t\int \frac{R\zeta-p_+\phi}{v}U_x(\zeta-\Theta)_+dxds
+\k\int_0^t\int_{\Omega_2}\frac{\zeta_x\Theta_x}{V}dxds
-\k\int_0^t\int_{\Omega_2}\frac{\phi\Theta_x^2}{vV}dxds
\\[2mm]
\di+2\mu\int_0^t\int\frac{\psi_xU_x}{V}(\zeta-\Theta)_+dxds
-\mu\int_0^t\int\frac{\phi U_ x^2}{vV}(\zeta-\Theta)_+dxds
+2\int_0^t\int_{\Omega_2}\frac{R\zeta-p_+\phi}{v}\psi\zeta_xdxds\\[3mm]
\di-2\int_0^t\int_{\Omega_2}\frac{R\zeta-p_+\phi}{v}\psi\Theta_x dxds
-2\mu\int_0^t\int_{\Omega_2}\frac{\psi\psi_x\zeta_x}{v}dxds
+2\mu\int_0^t\int_{\Omega_2}\frac{\phi\psi U_x}{vV}\zeta_xdxds\\[3mm]
\di+2\mu\int_0^t\int_{\Omega_2}\frac{\psi\psi_x\Theta_x}{v}dxds
-2\mu\int_0^t\int_{\Omega_2}\frac{\phi\psi}{vV}U_x\Theta_xdxds
-2\int_0^t\int\psi\widetilde R_1(\zeta-\Theta)_+dxds\\[3mm]
\di-\int_0^t\int\widetilde R_2(\zeta-\Theta)_+dxds
-c_{\nu}\int_0^t\int\partial_t\Theta(\zeta-\Theta)_+dxds
-\int_0^t\int_{\Omega_2}\psi^2\partial_t\Theta dxds\\[3mm]
\di+\frac{\mu}{c_{\nu}}\int_0^t\int_{\Omega_2}\psi^2\left(\frac{\psi_x^2+2\psi_xU_x}{v}-\frac{\phi U_x^2}{vV}\right)dxds
-\frac{1}{c_{\nu}}\int_0^t\int_{\Omega_2}\psi^2\widetilde R_2dxds\\[3mm]
\di-\frac{1}{c_{\nu}}\int_0^t\int_{\Omega_2}\psi^2\left (\frac{R\zeta+R\Theta}{v}\psi_x+\frac{R\zeta-p_+\phi}{v}U_x\right)dxds
+\frac{\k}{c_{\nu}}\int_0^t\int_{\Omega_2}\psi^2\left(\frac{\t_x}{v}-\frac{\Theta_x}{V}\right)_xdxds\\[3mm]
\di\triangleq \int\left(\frac{c_{\nu}}{2}(\zeta_0(x)-\Theta(x,0))_+^2+\psi_0^2(x)(\zeta_0(x)-\Theta(x,0))_+\right)dx
+\sum_{i=1}^{20}I_i.
\end{array}
\end{equation}
We will estimate \eqref{long} term by term. Recalling \eqref{basic}, \eqref{v} and \eqref{omega-bound}, it holds that
\begin{equation}
\begin{array}{ll}\label{i1}
\di|I_1|&\di=\left|\int_0^t\int\frac{R\zeta-2p_+\phi-R\Theta}{v}\psi_x(\zeta-\Theta)_+ dxds\right|\\[2mm]
&\di\leq \frac{\mu}{4}\int_0^t\int\frac{\psi_x^2}{v}(\zeta-\Theta)_+ dxds
+C_0\int_0^t\int(\zeta^2+\phi^2+1)(\zeta-\Theta)_+dxds\\[2mm]
&\di\leq \frac{\mu}{4}\int_0^t\int\frac{\psi_x^2}{v}(\zeta-\Theta)_+ dxds
+C_0\int_0^t\int(\zeta^2+\phi^2\zeta)(\zeta-\Theta)_+dxds\\
&\di\leq \frac{\mu}{4}\int_0^t\int\frac{\psi_x^2}{v}(\zeta-\Theta)_+ dxds
+C_0\int_0^t\int(\zeta+\phi^2)\left(\zeta-\frac{1}{2}\Theta\right)_+^2dxds\\[2mm]
&\di\leq \frac{\mu}{4}\int_0^t\int\frac{\psi_x^2}{v}(\zeta-\Theta)_+ dxds
+C_0\int_0^t\max_{x\in\mathbb{R}}\left(\zeta-\frac{1}{2}\Theta\right)_+^2\int_{\{\zeta>\frac{\Theta}{2}\}}(\zeta+\phi^2)dxds\\
&\di\leq \frac{\mu}{4}\int_0^t\int\frac{\psi_x^2}{v}(\zeta-\Theta)_+ dxds
+C_0\int_0^t\max_{x\in\mathbb{R}}\left(\zeta-\frac{1}{2}\Theta\right)_+^2ds.
\end{array}
\end{equation}
It follows from \eqref{important}, \eqref{basic}, \eqref{phi-x3} and Cauchy's inequality, one has
\begin{equation}
\begin{array}{ll}\label{i3}
\di|I_2|+|I_3|+|I_4|\\
\di\leq C_0\int_0^t\int(\phi^2+\zeta^2)|U_x|dxds+
\k\int_0^t\int_{\Omega_2}\left|\frac{\zeta}{\Theta}\right|\left(\left|\frac{\zeta_x\Theta_x}{V}\right|
+\left|\frac{\phi\Theta_x^2}{vV}\right|\right) dxds\\
\di\leq \frac{\k}{8}\int_0^t\int_{\Omega_2}\frac{\zeta_x^2}{v}dxds
+C_0\int_0^t\int(\phi^2+\zeta^2)(|U_x|+\Theta_x^2)dxds\\
\di\leq \frac{\k}{8}\int_0^t\int_{\Omega_2}\frac{\zeta_x^2}{v}dxds+C_0.
\end{array}
\end{equation}
Similarly,
\begin{equation}
\begin{array}{ll}\label{i4}
\di|I_5|+|I_6|+|I_{10}|+|I_{11}|+|I_{12}|\\
\di\leq C_0\int_0^t\int\frac{\psi_x^2}{\t}dxds
+\frac{\k}{8}\int_0^t\int_{\Omega_2}\frac{\zeta_x^2}{v}dxds\\
\di\quad+C(C_0,M)\int_0^t\int(\phi^2+\psi^2+\zeta^2)(U_x^2+\Theta_x^2)dxds\\
\di\leq\frac{\k}{8}\int_0^t\int_{\Omega_2}\frac{\zeta_x^2}{v}dxds+C_0.
\end{array}
\end{equation}
By Cauchy's inequality, \eqref{basic}, \eqref{v} and \eqref{omega-bound}, it holds that
\begin{equation}
\begin{array}{ll}\label{i7}
\di|I_7|\leq\frac{\k}{8}\int_0^t\int_{\Omega_2}\frac{\zeta_x^2}{v}dxds
+C_0\int_0^t\int_{\Omega_2}(\zeta^2+\phi^2)\psi^2dxds\\
\di\quad\leq\frac{\k}{8}\int_0^t\int_{\Omega_2}\frac{\zeta_x^2}{v}dxds
+C_0\int_0^t\int_{\Omega_2}(\zeta^2\psi^2+\phi^2\psi^4+\phi^2)dxds\\
\di\quad\leq\frac{\k}{8}\int_0^t\int_{\Omega_2}\frac{\zeta_x^2}{v}dxds
+C_0\int_0^t\int_{\Omega_2}(\zeta^2(\phi^2+\psi^2)+\phi^2\psi^4)dxds\\
\di\quad\leq\frac{\k}{8}\int_0^t\int_{\Omega_2}\frac{\zeta_x^2}{v}dxds
+C_0\int_0^t\int_{\Omega_2}\left(\left(\zeta-\frac{1}{2}\Theta\right)^2(\phi^2+\psi^2)+\phi^2\psi^4\right)dxds\\
\di\quad\leq\frac{\k}{8}\int_0^t\int_{\Omega_2}\frac{\zeta_x^2}{v}dxds
+C_0\int_0^t\left(\max_{x\in\mathbb{R}}\left(\zeta-\frac{1}{2}\Theta\right)_+^2+\max_{x\in\mathbb{R}}\psi^4\right)\int(\phi^2+\psi^2)dxds\\
\di\quad\leq\frac{\k}{8}\int_0^t\int_{\Omega_2}\frac{\zeta_x^2}{v}dxds
+C_0\int_0^t\left(\max_{x\in\mathbb{R}}\left(\zeta-\frac{1}{2}\Theta\right)_+^2+\max_{x\in\mathbb{R}}\psi^4\right)ds.
\end{array}
\end{equation}
Similarly, one has
\begin{equation}
\begin{array}{ll}\label{i8}
\di|I_8|\leq C_0\int_0^t\int_{\Omega_2}(\psi^2+(\phi^2+\zeta^2)\Theta_x^2)dxds\\
\di\quad\leq C_0\int_0^t\int_{\Omega_2}(\zeta^2\psi^2+(\phi^2+\zeta^2)\Theta_x^2)dxds\\
\di\quad\leq C_0\int_0^t\max_{x\in\mathbb{R}}\left(\zeta-\frac{1}{2}\Theta\right)_+^2\int\psi^2dxds
+C_0\int_0^t\int(\phi^2+\zeta^2)\Theta_x^2dxds\\
\di\quad\leq C_0\int_0^t\max_{x\in\mathbb{R}}\left(\zeta-\frac{1}{2}\Theta\right)_+^2ds+C_0.
\end{array}
\end{equation}
Using Cauchy's inequality, it holds that
\begin{equation}
\begin{array}{ll}\label{i9}
\di|I_9|\leq \frac{\k}{8}\int_0^t\int_{\Omega_2}\frac{\zeta_x^2}{v}dxds
+C_0\int_0^t\int\psi^2\psi_x^2dxds.
\end{array}
\end{equation}
Recalling Lemma \ref{decay}, \eqref{basic} and \eqref{phi-x3}, one has
\begin{equation}
\begin{array}{ll}\label{i10}
\di|I_{13}|+|I_{14}|+|I_{15}|+|I_{16}|+|I_{18}|\\
\di\leq C\d \int_0^t\int (1+s)^{-1}(\psi^2+\zeta^2)e^{-\frac{c_1x^2}{1+s}}dxds
\leq C_0.
\end{array}
\end{equation}
By Cauchy's inequality, and using \eqref{basic}, \eqref{phi-x3} and Lemma \ref{decay}, we obtain
\begin{equation}
\begin{array}{ll}\label{i17}
\di|I_{17}|\leq C_0\int_0^t\int(\psi^2\psi_x^2+\psi^2|\psi_x||U_x|+\psi^2|\phi|U_x^2)dxds\\
\di\quad\leq C_0\int_0^t\int(\psi^2\psi_x^2+\psi^2U_x^2+\psi^2|\phi|U_x^2)dxds\\
\di\quad\leq C_0\int_0^t\int\psi^2\psi_x^2dxds+C(C_0,M)\int_0^t\int\psi^2U_x^2dxds\\
\di\quad\leq C_0\int_0^t\int\psi^2\psi_x^2dxds+C_0.
\end{array}
\end{equation}
Similarly,
\begin{equation}
\begin{array}{ll}\label{i19}
\di|I_{19}|\leq\int_0^t\int_{\Omega_2}(\psi^2\psi_x^2+\psi^2\zeta^2+\psi^4|U_x|+(\zeta^2+\phi^2)|U_x|)dxds\\
\di\leq C_0\int_0^t\int\psi^2\psi_x^2dxds
+C_0\int_0^t\max_{x\in\mathbb{R}}\left(\zeta-\frac{1}{2}\Theta\right)_+^2\int\psi^2dxds\\
\di\quad+\int_0^t\max_{x\in\mathbb{R}}\psi^4\int|U_x|dxds+C_0\\
\di\leq C_0\int_0^t\int\psi^2\psi_x^2dxds
+C_0\int_0^t\left(\max_{x\in\mathbb{R}}\left(\zeta-\frac{1}{2}\Theta\right)_+^2+\max_{x\in\mathbb{R}}\psi^4\right)ds
+C_0.
\end{array}
\end{equation}
Finally, for
\begin{equation}\label{cut-off}
\di \varphi_{\eta}(z)=\left\{
\begin{array}{ll}
1, & z>\eta,\\
z/\eta, & 0< z\leq \eta,\\
0  & z\leq 0.
\end{array}
\right.
\end{equation}
Integrating by parts shows
\begin{equation}
\begin{array}{ll}\label{i20}
\di I_{20}=\frac{\k}{c_{\nu}}\int_0^t\int_{\Omega_2}\psi^2\left(\frac{\t_x}{v}-\frac{\Theta_x}{V}\right)_xdxds\\
\di=\frac{\k}{c_{\nu}}\lim_{\eta\rightarrow0+}\int_0^t\int\varphi_{\eta}(\zeta-\Theta)\psi^2\left(\frac{\zeta_x-\Theta_x}{v}\right)_xdxds\\
\di\quad+\frac{\k}{c_{\nu}}\int_0^t\int_{\Omega_2}\psi^2\left(\frac{2\Theta_x}{v}-\frac{\Theta_x}{V}\right)_xdxds
\triangleq I_{20}^1+I_{20}^2.
\end{array}
\end{equation}
We have
\begin{equation}
\begin{array}{ll}\label{i201}
\di I_{20}^1&\di=-\frac{2\k}{c_{\nu}}\lim_{\eta\rightarrow 0+}\int_0^t\int\varphi_{\eta}(\zeta-\Theta)\psi\psi_x\frac{\zeta_x-\Theta_x}{v}dxds\\[3mm]
&\di\quad-\frac{\k}{c_{\nu}}\lim_{\eta\rightarrow 0+}\int_0^t\int\varphi_{\eta}'(\zeta-\Theta)\frac{\psi^2(\zeta_x-\Theta_x)^2}{v}dxds\\[3mm]
&\di\leq C\int_0^t\int_{\Omega_2}\frac{|\psi\psi_x\zeta_x|+|\psi\psi_x\Theta_x|}{v} dxds\\[3mm]
&\di\leq\frac{\k}{8}\int_0^t\int_{\Omega_2}\frac{\zeta_x^2}{v}dxds+C_0\int_0^t\int\psi^2\psi_x^2dxds\\[3mm]
&\di\quad+\int_0^t\int\frac{\psi_x^2}{\t v}dxds+C(C_0,M)\int_0^t\int\psi^2\Theta_x^2dxds\\[3mm]
&\di\leq\frac{\k}{8}\int_0^t\int_{\Omega_2}\frac{\zeta_x^2}{v}dxds+C_0\int_0^t\int\psi^2\psi_x^2dxds+C_0,
\end{array}
\end{equation}
where in the second inequality we have used both $\varphi_{\eta}(z)\in [0,1]$ and $\varphi_{\eta}'(z)\geq 0$.
Similarly,
\begin{equation}
\begin{array}{ll}\label{i202}
\di I_{20}^2=\frac{\k}{c_{\nu}}\int_0^t\int_{\Omega_2}\psi^2\left(\frac{2\Theta_{xx}}{v}-\frac{\Theta_{xx}}{V}
-\frac{2\Theta_xv_x}{v^2}+\frac{\Theta_xV_x}{V^2}\right) dxds\\
\di\quad\leq C_0\int_0^t\int\psi^2(|\Theta_{xx}|+\Theta_x^2)+\psi^2|\phi_x||\Theta_x| dxds\\
\di\quad\leq C_0+C_0\int_0^t\max_{x\in\mathbb{R}}\psi^4\left(\int|\Theta_x|dx\right)ds
+C(C_0,M)\d\int_0^t\int \frac{\theta\phi_x^2}{v^3}dxds\\
\di\quad\leq C_0+C_0 \int_0^t\max_{x\in\mathbb{R}}\psi^4ds.
\end{array}
\end{equation}
Noticing that
\begin{equation}
\begin{array}{ll}\label{5.4}
\di\int_0^t\int(\t\psi_x^2+\zeta_x^2)dxds=\int_0^t\int_{\{\zeta>2\Theta\}}(\t\psi_x^2+\zeta_x^2)dxds
+\int_0^t\int_{\{\zeta\leq 2\Theta\}}(\t\psi_x^2+\zeta_x^2)dxds\\[3mm]
\di\quad\leq\int_0^t\int_{\{\zeta>2\Theta\}}\left(\frac{3}{2}\psi_x^2\zeta+C_0\frac{\zeta_x^2}{v}\right) dxds
+\int_0^t\int_{\{\zeta\leq 2\Theta\}}\left(\frac{\psi_x^2}{\t}+\frac{\zeta_x^2}{\t^2}\right)\t^2dxds\\[3mm]
\di\quad\leq\int_0^t\int_{\{\zeta>2\Theta\}}\left(3\psi_x^2(\zeta-\Theta)+C_0\frac{\zeta_x^2}{v}\right) dxds
+C\int_0^t\int_{\{\zeta\leq 2\Theta\}}\left(\frac{\psi_x^2}{\t}+\frac{\zeta_x^2}{\t^2}\right)dxds\\[3mm]
\di\quad\leq C_0\int_0^t\int_{\Omega_2}\left(\frac{\psi_x^2}{v}(\zeta-\Theta)_+ +\frac{\zeta_x^2}{v}\right)dxds+C_0.
\end{array}
\end{equation}
Substituting the estimates \eqref{i1}-\eqref{i202} into \eqref{long}, and using \eqref{5.4}, we have
\begin{equation}
\begin{array}{ll}\label{5.3}
\di \int(\zeta-\Theta)_+^2 dx+\int_0^t\int(\t\psi_x^2+\zeta_x^2)dxds\leq C_0\\
\di\quad +C_0\int_0^t\left(\max_{x\in\mathbb{R}}\left(\zeta-\frac{1}{2}\Theta\right)_+^2+\max_{x\in\mathbb{R}}\psi^4\right)ds
+C_0\int_0^t\int\psi^2\psi_x^2dxds.
\end{array}
\end{equation}

\

\underline{ Step 2. }\quad To estimate the last term on the right hand side of \eqref{5.3},
we multiply $\eqref{perturb}_2$ by $\psi^3$, and integrate the resulted equation over $\mathbb{R}\times[0,t]$ to get
\begin{equation}
\begin{array}{ll}\label{5.5}
\di\frac{1}{4}\int\psi^4dx+3\mu\int_0^t\int\frac{\psi^2\psi_x^2}{v} dxds
=\frac{1}{4}\int\psi_0^4dx+3R\int_0^t\int\frac{\zeta\psi^2\psi_x}{v}dxds\\[3mm]
\di\quad-3p_+\int_0^t\int\frac{\phi\psi^2\psi_x}{v}dxds+3\mu\int_0^t\int\frac{\phi U_x}{vV}\psi^2\psi_x dxds
-\int_0^t\int\widetilde R_1\psi^3dxds\\
\di\quad\triangleq\frac{1}{4}\int\psi_0^4dx+\sum_{i=1}^4J_i.
\end{array}
\end{equation}
It follows from \eqref{basic} and \eqref{v} that,
\begin{equation}
\begin{array}{ll}\label{j1}
\di|J_1|&\di=3R\int_0^t\int_{\{\zeta>\Theta\}}\frac{\zeta\psi^2\psi_x}{v}dxds
+3R\int_0^t\int_{\{\zeta\leq\Theta\}}\frac{\zeta\psi^2\psi_x}{v}dxds\\[3mm]
&\di\leq \mu\int_0^t\int_{\{\zeta>\Theta\}}\frac{\psi^2\psi_x^2}{v}dxds+C_0\int_0^t\int_{\{\zeta>\Theta\}}\zeta^2\psi^2dxds\\[3mm]
&\di\quad+\int_0^t\int_{\{\zeta\leq\Theta\}}\psi_x^2dxds+C_0\int_0^t\int_{\{\zeta\leq\Theta\}}\zeta^2\psi^4dxds\\[3mm]
&\di\leq\mu\int_0^t\int\frac{\psi^2\psi_x^2}{v}dxds+C_0\int_0^t\max_{x\in\mathbb{R}}\left(\zeta-\frac{1}{2}\Theta\right)_+^2\left(\int\psi^2dx\right)ds\\[3mm]
&\di\quad+C\int_0^t\int_{\{\zeta\leq\Theta\}}\frac{\psi_x^2}{\t}dxds
+C_0\int_0^t\max_{x\in\mathbb{R}}\psi^4\left(\int_{\{\zeta\leq\Theta\}}\zeta^2dx\right)ds\\[3mm]
&\di\leq\mu\int_0^t\int\frac{\psi^2\psi_x^2}{v}dxds
+C_0\int_0^t\left(\max_{x\in\mathbb{R}}\left(\zeta-\frac{1}{2}\Theta\right)_+^2+\max_{x\in\mathbb{R}}\psi^4\right)ds+C_0.
\end{array}
\end{equation}
Recalling \eqref{basic}, \eqref{v}, and using Cauchy's inequality, it holds that
\begin{equation}
\begin{array}{ll}\label{j2}
\di|J_2|&\di\leq \v\int_0^t\int\psi_x^2dxds+C(\v^{-1},C_0)\int_0^t\int\phi^2\psi^4dxds\\
&\di\leq\v\int_0^t\int\psi_x^2dxds+C(\v^{-1},C_0)\int_0^t\max_{x\in\mathbb{R}}\psi^4\left(\int\phi^2dx\right)ds\\
&\di\leq\v\int_0^t\int\psi_x^2dxds+C(\v^{-1},C_0)\int_0^t\max_{x\in\mathbb{R}}\psi^4ds.
\end{array}
\end{equation}
From \eqref{basic}, \eqref{important} and \eqref{phi-x3}, one has
\begin{equation}
\begin{array}{ll}\label{j3}
\di|J_3|&\di\leq \mu\int_0^t\int\frac{\psi^2\psi_x^2}{v}dxds
+C_0\int_0^t\int\phi^2\psi^2U_x^2dxds\\
&\di\leq\mu\int_0^t\int\frac{\psi^2\psi_x^2}{v}dxds
+C(C_0,M)\int_0^t\int\phi^2U_x^2dxds\\
&\di\leq\mu\int_0^t\int\frac{\psi^2\psi_x^2}{v}dxds+C_0.
\end{array}
\end{equation}
By Lemma \ref{decay}, \eqref{basic}, \eqref{important} and \eqref{phi-x3}, we have
\begin{equation}
\begin{array}{ll}\label{j4}
\di|J_4|&\di\leq O(1)\d\int_0^t\int(1+s)^{-1}|\psi|^3e^{-\frac{c_1x^2}{1+s}}dxds\\
&\di\leq C(M)\d\int_0^t\int(1+s)^{-1}\psi^2e^{-\frac{c_1x^2}{1+s}}dxds\leq C_0.
\end{array}
\end{equation}
Putting the estimates \eqref{j1}-\eqref{j4} into \eqref{5.5} gives
\begin{equation}
\begin{array}{ll}\label{5.6}
\di\int\psi^4dx+\int_0^t\int\psi^2\psi_x^2dxds\leq C_0+C_0\v\int_0^t\int\psi_x^2dxds\\
\di\qquad +C(\v^{-1},C_0)\int_0^t\left(\max_{x\in\mathbb{R}}\left(\zeta-\frac{1}{2}\Theta\right)_+^2+\max_{x\in\mathbb{R}}\psi^4\right)ds.
\end{array}
\end{equation}
Noticing that
\begin{equation}\label{psi-x}
\di 2\int_0^t\int\psi_x^2dxds\leq \int_0^t\int\frac{\psi_x^2}{\t}dxds+\int_0^t\int\t\psi_x^2dxds
\leq C_0+C_0\int_0^t\int\t\psi_x^2dxds.
\end{equation}
Combining \eqref{5.3} and \eqref{5.6}, choosing $\v$ suitable small, we have
\begin{equation}
\begin{array}{ll}\label{5.7}
\di \int((\zeta-\Theta)_+^2+\psi^4) dx+\int_0^t\int((\psi^2+\t)\psi_x^2+\zeta_x^2)dxds\\
\di\quad \leq C_0+C_0\int_0^t\left(\max_{x\in\mathbb{R}}\left(\zeta-\frac{1}{2}\Theta\right)_+^2+\max_{x\in\mathbb{R}}\psi^4\right)ds.
\end{array}
\end{equation}

\

\underline{ Step 3. }\quad It remains to estimate the last two terms on the right hand side of \eqref{5.7}.
For $x\in\mathbb{R}$,
\begin{equation}
\begin{array}{ll}
\di\left(\zeta-\frac{1}{2}\Theta\right)_+^2=\int_{-\infty}^x2\left(\zeta-\frac{1}{2}\Theta\right)_+\left(\zeta_x-\frac{1}{2}\Theta_x\right)dx\\[3mm]
\di\quad\leq C\int\left(\zeta-\frac{1}{2}\Theta\right)_+(|\zeta_x|+|\Theta_x|)dx\\[3mm]
\di\quad\leq \v\int\left(\zeta-\frac{1}{2}\Theta\right)_+^2\t dx
+\frac{C}{\v}\int_{\{\zeta>\frac{\Theta}{2}\}}\left(\frac{\zeta_x^2}{\t}+\frac{\Theta_x^2}{\t}\right)dx\\[4mm]
\di\quad\leq 3\v\int\left(\zeta-\frac{1}{2}\Theta\right)_+^2\zeta dx
+\frac{C}{\v}\int\frac{\zeta_x^2}{\t}dx+\frac{C}{\v}\int_{\{\zeta>\frac{\Theta}{2}\}}\zeta^2\Theta_x^2dx\\[3mm]
\di\quad\leq 3\v\max_{x\in\mathbb{R}}\left(\zeta-\frac{1}{2}\Theta\right)_+^2\int_{\{\zeta>\frac{\Theta}{2}\}}\zeta dx
+\frac{C}{\v}\int\frac{\zeta_x^2}{\t}dx+\frac{C}{\v}\int\zeta^2\Theta_x^2dx\\[3mm]
\di\quad\leq \v C_0\max_{x\in\mathbb{R}}\left(\zeta-\frac{1}{2}\Theta\right)_+^2+\frac{C}{\v}\int\frac{\zeta_x^2}{\t}dx
+\frac{C}{\v}\int\zeta^2\Theta_x^2dx.
\end{array}
\end{equation}
This yields
\begin{equation}\label{max-theta}
\di\max_{x\in\mathbb{R}}\left(\zeta-\frac{1}{2}\Theta\right)_+^2\leq C_0\int\frac{\zeta_x^2}{\t}dx
+C_0\int\zeta^2\Theta_x^2dx.
\end{equation}

\begin{equation}
\begin{array}{ll}
\di\psi^4=\int_{-\infty}^x4\psi^3\psi_xdx\leq4\int_{\{\zeta>\Theta\}}|\psi|^3|\psi_x|dx
+4\int_{\{\zeta\leq\Theta\}}|\psi|^3|\psi_x|dx\\[3mm]
\di\quad\leq\v\int_{\{\zeta>\Theta\}}|\psi|^5\sqrt{\t}dx+\frac{C}{\v}\int_{\{\zeta>\Theta\}}\psi_x^2\frac{|\psi|}{\sqrt{\t}}dx
+\v\int_{\{\zeta\leq\Theta\}}\psi^6\t dx+\frac{C}{\v}\int_{\{\zeta\leq\Theta\}}\frac{\psi_x^2}{\t}dx\\[3mm]
\di\quad\leq\v\max_{x\in\mathbb{R}}\psi^4\int_{\{\zeta>\Theta\}}(\psi^2+\theta)dx
+C\v\max_{x\in\mathbb{R}}\psi^4\int_{\{\zeta\leq\Theta\}}\psi^2dx
+\frac{C}{\v}\int\psi_x^2\left(\frac{|\psi|}{\sqrt{\t}}+\frac{1}{\t}\right)dx\\[3mm]
\di\quad\leq \v C_0\max_{x\in\mathbb{R}}\psi^4+\frac{C}{\v}\int\psi_x^2\left(\frac{|\psi|}{\sqrt{\t}}+\frac{1}{\t}\right)dx,
\end{array}
\end{equation}
which directly gives
\begin{equation}\label{max-psi}
\di\max_{x\in\mathbb{R}}\psi^4\leq C_0\int\psi_x^2
\left(\frac{|\psi|}{\sqrt{\t}}+\frac{1}{\t}\right)dx.
\end{equation}
Substituting \eqref{max-theta} and \eqref{max-psi} into \eqref{5.7},
recalling \eqref{basic}, \eqref{important} and \eqref{phi-x3}, and choosing $\v$ suitable small, it holds that
\begin{equation}
\begin{array}{ll}\label{5.8}
\di \sup_{0\leq t\leq T}\int((\zeta-\Theta)_+^2+\psi^4) dx+\int_0^T\int((\psi^2+\t)\psi_x^2+\zeta_x^2)dxdt\\
\di\quad \leq C_0+C_0\int_0^T\int\left(\frac{\zeta_x^2}{\t}+\psi_x^2\left(\frac{|\psi|}{\sqrt{\t}}+\frac{1}{\t}\right)\right)dxdt\\
\di\quad\leq C_0+\frac{1}{2}\int_0^T\int(\zeta_x^2+\psi^2\psi_x^2)dxdt
+C_0\int_0^T\int\left(\frac{\zeta_x^2}{\t^2}+\frac{\psi_x^2}{\t}\right)dxdt\\
\di\quad\leq C_0+\frac{1}{2}\int_0^T\int(\zeta_x^2+\psi^2\psi_x^2)dxdt.
\end{array}
\end{equation}
From\eqref{basic}, \eqref{omega-bound}, we have
\begin{equation}\label{5.9}
\di\int_{\{\zeta\leq2\Theta\}}\zeta^2dx\leq C\int_{\mathbb{R}}\Phi\left(\frac{\t}{\Theta}\right)dx\leq C_0,
\end{equation}
and
\begin{equation}\label{5.10}
\di\int_{\{\zeta>2\Theta\}}\zeta^2dx\leq4\int_{\{\zeta>2\Theta\}}(\zeta-\Theta)^2dx
\leq 4\int(\zeta-\Theta)_+^2dx.
\end{equation}
Thus combining \eqref{5.8}-\eqref{5.10}, the proof of Lemma \ref{lemma-ll-1} is completed.
\hfill $\Box$

\

\begin{lemma}\label{high-order}
Suppose that $(\phi,\psi,\zeta)\in X([0,T])$ satisfies $\d=|\t_+-\t_-|\leq\d_0$
with suitable small $\d_0$, it holds
\begin{equation}\label{high-deri}
\di\sup_{0\leq t\leq T}\int(\phi_x^2+\psi_x^2+\zeta_x^2)dx
+\int_0^T\int(\t\phi_x^2+\psi_{xx}^2+\zeta_{xx}^2)dxdt\leq C_0.
\end{equation}
\end{lemma}
\textbf{Proof}:\quad Due to \eqref{basic}, \eqref{v} and \eqref{new-nergy}, some terms of \eqref{phi-x1} can be considered more carefully, that is,
\begin{equation}
\begin{array}{ll}
\di \left|\frac{R}{v}\zeta_x\frac{\tilde v_x}{\tilde v}\right|
+\left|\psi_xU_x\left(\frac{1}{v}-\frac{1}{V}\right)\right|+\frac{\psi_x^2}{v}\\[3mm]
\di \leq \frac{R\t}{4v}\left(\frac{\tilde v_x}{\tilde v}\right)^2+
C_0\frac{\zeta_x^2}{\t}+C_0\psi_x^2+C_0\phi^2U_x^2\\[3mm]
\di \leq \frac{R\t}{4v}\left(\frac{\tilde v_x}{\tilde v}\right)^2+
C_0\left(\zeta_x^2+\frac{\zeta_x^2}{\t^2}+\t\psi_x^2+\frac{\psi_x^2}{\t}\right)
+C_0\phi^2U_x^2.
\end{array}
\end{equation}
The other terms in \eqref{phi-x1} can be estimated the same as in step 2 in Lemma \ref{basic-lemma}.
Integrating \eqref{phi-x1} over $\mathbb{R}\times(0,t)$, recalling \eqref{basic} and \eqref{new-nergy}, we have
\begin{equation}\label{new-phi-x}
\di \sup_{0\leq t\leq T}\int\phi_x^2 dx+\int_0^T\int\t\phi_x^2dxdt\leq C_0.
\end{equation}
Multiplying $\eqref{perturb}_2$ by $-\psi_{xx}$, integrating the resulted equation over $\mathbb{R}\times(0,t)$,
and noticing that
\begin{equation*}
\begin{array}{ll}
\di(p-p_+)_x&\di=\left(\frac{R\zeta-p_+\phi}{v}\right)_x=\frac{R\zeta_x-p_+\phi_x}{v}
-\frac{R\zeta-p_+\phi}{v^2}\phi_x-\frac{R\zeta-p_+\phi}{v^2}V_x\\[3mm]
&\di=\frac{R\zeta_x}{v}-\frac{R\t\phi_x}{v^2}-\frac{R\zeta-p_+\phi}{v^2}V_x.
\end{array}
\end{equation*}
Then we have
\begin{equation}
\begin{array}{ll}\label{6.1}
\di\int\frac{\psi_x^2}{2}dx+\mu\int_0^t\int\frac{\psi_{xx}^2}{v}dxds
=\int\frac{\psi_{0x}^2}{2}dx+\int_0^t\int\left(\frac{R\zeta_x}{v}-\frac{R\t\phi_x}{v^2}
-\frac{R\zeta-p_+\phi}{v^2}V_x\right)\psi_{xx}dxds\\
\di\quad-\mu\int_0^t\int\psi_x\left(\frac{1}{v}\right)_x\psi_{xx}dxds
+\mu\int_0^t\int\left(\frac{U_x}{V}-\frac{U_x}{v}\right)_x\psi_{xx} dxds
+\int_0^t\int\widetilde R_1\psi_{xx}dxds.
\end{array}
\end{equation}
In the following, each term on the right hand side of \eqref{6.1} will be estimated.
From \eqref{important}, \eqref{new-nergy} and \eqref{new-phi-x}, one has
\begin{equation}
\begin{array}{ll}\label{psi-xx-1}
\di\left|\int_0^t\int\left(\frac{R\zeta_x}{v}-\frac{R\t\phi_x}{v^2}
-\frac{R\zeta-p_+\phi}{v^2}V_x\right)\psi_{xx}dxds\right|\\
\di\leq \frac{\mu}{8}\int_0^t\int\frac{\psi_{xx}^2}{v}dxds
+C_0\int_0^t\int\zeta_x^2+\theta^2\phi_x^2+(\phi^2+\zeta^2)\Theta_x^2dxds\\
\di\leq \frac{\mu}{8}\int_0^t\int\frac{\psi_{xx}^2}{v}dxds
+C_0+\max_{x,t}\t\int_0^t\int\t\phi_x^2dxds\\
\di\leq \frac{\mu}{8}\int_0^t\int\frac{\psi_{xx}^2}{v}dxds
+C_0+C_0\max_{x,t}\t.
\end{array}
\end{equation}
By Cauchy's inequality and Sobolev's inequality, and recalling \eqref{basic}, \eqref{new-nergy},
\eqref{psi-x} and \eqref{new-phi-x}, we obtain
\begin{equation}
\begin{array}{ll}\label{psi-xx-2}
\di\left|\mu\int_0^t\int\psi_x\left(\frac{1}{v}\right)_x\psi_{xx}dxds\right|
\leq C\int_0^t\int\left|\frac{\psi_x\phi_x\psi_{xx}}{v^2}\right|+\left|\frac{\psi_xV_x\psi_{xx}}{v^2}\right|dxds\\
\di\leq\frac{\mu}{8}\int_0^t\int\frac{\psi_{xx}^2}{v}dxds
+C_0\int_0^t\int \psi_x^2\phi_x^2+\psi_x^2\Theta_x^2dxds\\
\di\leq\frac{\mu}{8}\int_0^t\int\frac{\psi_{xx}^2}{v}dxds
+C_0\int_0^t\|\psi_x\|^2_{L^{\infty}}\|\phi_x\|^2ds
+C_0\int_0^t\int\psi_x^2dxds\\
\di\leq\frac{\mu}{8}\int_0^t\int\frac{\psi_{xx}^2}{v}dxds
+C_0\int_0^t\|\psi_x\|\|\psi_{xx}\|ds
+C_0\int_0^t\int\psi_x^2dxds\\
\di\leq\frac{\mu}{4}\int_0^t\int\frac{\psi_{xx}^2}{v}dxds
+C_0\int_0^t\int\psi_x^2dxds\\
\di\leq\frac{\mu}{4}\int_0^t\int\frac{\psi_{xx}^2}{v}dxds+C_0.
\end{array}
\end{equation}
Similarly,
\begin{equation}
\begin{array}{ll}\label{psi-xx-3}
\di\mu\int_0^t\int\left(\frac{U_x}{V}-\frac{U_x}{v}\right)_x\psi_{xx} dxds\\
\di=\mu\int_0^t\int\left(\frac{U_{xx}}{V}-\frac{U_{xx}}{v}-\frac{U_xV_x}{V^2}+\frac{v_xU_x}{v^2}\right)\psi_{xx}dxds\\
\di=\mu\int_0^t\int\left(\frac{\phi U_{xx}}{vV}-U_xV_x\frac{\phi(\phi+2V)}{v^2V^2}+\frac{\phi_xU_x}{v^2}\right)\psi_{xx}dxds\\
\di\leq\frac{\mu}{8}\int_0^t\int\frac{\psi_{xx}^2}{v}dxds
+C_0\int_0^t\int(\phi^2(U_{xx}^2+\Theta_x^2U_x^2)+\phi^4U_x^2\Theta_x^2+\phi_x^2U_x^2)dxds\\
\di\leq\frac{\mu}{8}\int_0^t\int\frac{\psi_{xx}^2}{v}dxds
+C_0+C(C_0,M)\d^2\int_0^t\int\t\phi_x^2dxds\\
\di\leq\frac{\mu}{8}\int_0^t\int\frac{\psi_{xx}^2}{v}dxds+C_0,
\end{array}
\end{equation}
and
\begin{equation}
\begin{array}{ll}\label{psi-xx-4}
\di\left|\int_0^t\int\widetilde R_1\psi_{xx}dxds\right|
\leq\frac{\mu}{8}\int_0^t\int\frac{\psi_{xx}^2}{v}dxds+C_0\int_0^t\int\widetilde R_1^2dxds\\
\di\leq\frac{\mu}{8}\int_0^t\int\frac{\psi_{xx}^2}{v}dxds+C_0\d^2\int_0^t\int(1+s)^{-3}e^{-\frac{c_1x^2}{1+s}}dxds\\
\di\leq\frac{\mu}{8}\int_0^t\int\frac{\psi_{xx}^2}{v}dxds+C_0.
\end{array}
\end{equation}
Substituting \eqref{psi-xx-1}-\eqref{psi-xx-4} into \eqref{6.1} shows
\begin{equation}\label{psi-xx}
\di\sup_{0\leq t\leq T}\int\psi_x^2dx+\int_0^T\int\psi_{xx}^2dxdt
\leq C_0+C_0\max_{x,t}\t.
\end{equation}

Multiplying $\eqref{perturb}_3$ by $-\zeta_{xx}$, then integrating the resulted equation over
$\mathbb{R}\times(0,t)$, we have
\begin{equation}
\begin{array}{ll}\label{6.2}
\di\frac{c_{\nu}}{2}\int\zeta_x^2dx+\k\int_0^t\int\frac{\zeta_{xx}^2}{v}dxds
=\frac{c_{\nu}}{2}\int\zeta_{0x}^2dx+\int_0^t\int(pu_x-p_+U_x)\zeta_{xx}dxds\\[3mm]
\di\quad-\k\int_0^t\int\zeta_x\left(\frac{1}{v}\right)_x\zeta_{xx}dxds
-\k\int_0^t\int\left(\frac{\Theta_x}{v}-\frac{\Theta_x}{V}\right)_x\zeta_{xx}dxds\\[3mm]
\di\quad-\mu\int_0^t\int\left(\frac{u_x^2}{v}-\frac{U_x^2}{V}\right)\zeta_{xx}dxds
+\int_0^t\int\widetilde R_2\zeta_{xx}dxds.
\end{array}
\end{equation}
We will estimate \eqref{6.2} one by one. First, we have
\begin{equation}
\begin{array}{ll}\label{zeta-xx-1}
\di\int_0^t\int(pu_x-p_+U_x)\zeta_{xx}dxds=
\int_0^t\int\left(\frac{R\t}{v}\psi_x+\frac{R\zeta-p_+\phi}{v}U_x\right)\zeta_{xx}dxds\\
\di\quad\leq\frac{\k}{8}\int_0^t\int\frac{\zeta_{xx}^2}{v}dxds+
C_0\int_0^t\int(\theta^2\psi_x^2+(\phi^2+\zeta^2)U_x^2)dxds\\
\di\quad\leq\frac{\k}{8}\int_0^t\int\frac{\zeta_{xx}^2}{v}dxds+
C_0\max_{x,t}\t\int_0^t\int\theta\psi_x^2dxds+
C_0\int_0^t\int(\phi^2+\zeta^2)U_x^2dxds\\
\di\quad\leq\frac{\k}{8}\int_0^t\int\frac{\zeta_{xx}^2}{v}dxds+
C_0\max_{x,t}\t+C_0.
\end{array}
\end{equation}
It follows from Cauchy's inequality, \eqref{new-nergy} and \eqref{new-phi-x} that
\begin{equation}
\begin{array}{ll}\label{zeta-xx-2}
\di-\k\int_0^t\int\zeta_x\left(\frac{1}{v}\right)_x\zeta_{xx}dxds
\leq C\int_0^t\int\frac{|\zeta_x\phi_x\zeta_{xx}|+|\zeta_xV_x\zeta_{xx}|}{v^2} dxds\\[3mm]
\di\quad\leq\frac{\k}{8}\int_0^t\int\frac{\zeta_{xx}^2}{v}dxds+
C_0\int_0^t\int\zeta_x^2\phi_x^2+\zeta_x^2\Theta_x^2dxds\\[3mm]
\di\quad\leq\frac{\k}{8}\int_0^t\int\frac{\zeta_{xx}^2}{v}dxds+
C_0\sup_{0\leq s\leq t}\|\phi_x\|^2\int_0^t\|\zeta_{x}\|\|\zeta_{xx}\|ds+
C_0\int_0^t\int\zeta_x^2dxds\\[3mm]
\di\quad\leq\frac{\k}{8}\int_0^t\int\frac{\zeta_{xx}^2}{v}dxds+
C_0\int_0^t\|\zeta_{x}\|\|\zeta_{xx}\|ds+C_0\\[3mm]
\di\quad\leq\frac{\k}{4}\int_0^t\int\frac{\zeta_{xx}^2}{v}dxds+C_0.
\end{array}
\end{equation}
Recalling Lemma \ref{decay} and \eqref{new-nergy}, and choosing $\delta$ suitable small, we have
\begin{equation}
\begin{array}{ll}\label{zeta-xx-3}
\di
-\k\int_0^t\int\left(\frac{\Theta_x}{v}-\frac{\Theta_x}{V}\right)_x\zeta_{xx}dxds\\[3mm]
\di=-\k\int_0^t\int\left(\frac{\Theta_{xx}}{v}-\frac{\Theta_{xx}}{V}
-\frac{\Theta_{x}v_x}{v^2}+\frac{\Theta_{x}V_x}{V^2}\right)\zeta_{xx}dxds\\[3mm]
\di=-\k\int_0^t\int\left(\frac{-\phi\Theta_{xx}}{vV}-\frac{\phi_x\Theta_x}{v^2}
+\frac{\phi(\phi+2V)}{v^2V^2}\Theta_xV_x\right)\zeta_{xx}dxds\\[3mm]
\di\leq\frac{\k}{8}\int_0^t\int\frac{\zeta_{xx}^2}{v}dxds+
C_0\int_0^t\int(\phi^2(\Theta_{xx}^2+\Theta_x^4)+\phi_x^2\Theta_x^2+\phi^4\Theta_x^4)dxds\\[3mm]
\di\leq\frac{\k}{8}\int_0^t\int\frac{\zeta_{xx}^2}{v}dxds+
C_0+C(C_0,M)\d^2\int_0^t\int\t\phi_x^2dxds+C(C_0,M)\int_0^t\int\phi^2\Theta_x^4dxds\\[3mm]
\di\leq\frac{\k}{8}\int_0^t\int\frac{\zeta_{xx}^2}{v}dxds+C_0.
\end{array}
\end{equation}
Similarly,
\begin{equation}
\begin{array}{ll}\label{zeta-xx-4}
\di-\mu\int_0^t\int\left(\frac{u_x^2}{v}-\frac{U_x^2}{V}\right)\zeta_{xx}dxds=
-\mu\int_0^t\int\left(\frac{\psi_x^2+2\psi_xU_x}{v}-\frac{\phi U_x^2}{vV}\right)\zeta_{xx}dxds\\[3mm]
\di\leq\frac{\k}{8}\int_0^t\int\frac{\zeta_{xx}^2}{v}dxds+
C_0\int_0^t\int(\psi_x^4+\psi_x^2U_x^2+\phi^2U_x^4)dxds\\[3mm]
\di\leq\frac{\k}{8}\int_0^t\int\frac{\zeta_{xx}^2}{v}dxds+
C_0\int_0^t\|\psi_x\|^3\|\psi_{xx}\|ds+C_0\\[3mm]
\di\leq\frac{\k}{8}\int_0^t\int\frac{\zeta_{xx}^2}{v}dxds+
\int_0^t\int\psi_{xx}^2dxds+C_0\sup_{0\leq s\leq t}\|\psi_x\|^4\int_0^t\int\psi_x^2dxds+C_0\\[3mm]
\di\leq\frac{\k}{8}\int_0^t\int\frac{\zeta_{xx}^2}{v}dxds+
C_0\max_{x,t}\t^2+C_0,
\end{array}
\end{equation}
and
\begin{equation}
\begin{array}{ll}\label{zeta-xx-5}
\di
\int_0^t\int\widetilde R_2\zeta_{xx}dxds\leq\frac{\k}{8}\int_0^t\int\frac{\zeta_{xx}^2}{v}dxds
+C_0\int_0^t\int\widetilde R_2^2dxds\\
\di\leq\frac{\k}{8}\int_0^t\int\frac{\zeta_{xx}^2}{v}dxds
+C_0\d^2\int_0^t\int(1+s)^{-4}e^{-\frac{c_1x^2}{1+s}}dxds\\
\di\leq\frac{\k}{8}\int_0^t\int\frac{\zeta_{xx}^2}{v}dxds+C_0.
\end{array}
\end{equation}
Substituting estimates \eqref{zeta-xx-1}-\eqref{zeta-xx-5} into \eqref{6.2} shows
\begin{equation}\label{zeta-xx}
\di\sup_{0\leq t\leq T}\int\zeta_x^2dx+\int_0^T\int\zeta_{xx}^2dxdt
\leq C_0+C_0\max_{x,t}\t^2.
\end{equation}

By Sobolev's inequality and \eqref{new-nergy}, \eqref{zeta-xx}, we have
\begin{equation}
\di \|\zeta\|^2_{L^{\infty}}\leq C\|\zeta\|\|\zeta_x\|\leq C_0+C_0\max_{x,t}\t.
\end{equation}
Noticing that
\begin{equation}
\di\max_{x,t}\t^2\leq 2\max_{x,t}\zeta^2+2\max_{x,t}\Theta^2\leq C_0+C_0\max_{x,t}\t.
\end{equation}
This yields
\begin{equation}\label{theta-upper}
\max_{x,t}\t\leq C_0,
\end{equation}
which, together with \eqref{new-phi-x}, \eqref{psi-xx} and \eqref{zeta-xx},
completes the proof of Lemma \ref{high-order}.
\hfill $\Box$

\

Finally, it follows from \eqref{new-nergy}, \eqref{high-deri} and equation $\eqref{perturb}_3$ that
$$
\di\int_0^{+\infty}\|\zeta_x\|^2+\left|\frac{d}{dt}\|\zeta_x\|^2\right| dt\leq C_0,
$$
which, together with  the Sobolev's inequality gives
\begin{equation}
\di\lim_{t\rightarrow \infty}\|\zeta\|^2_{L^{\infty}}
\leq C\lim_{t\rightarrow \infty}\|\zeta\|\|\zeta_x\|
\leq C_0\lim_{t\rightarrow \infty}\|\zeta_x\|=0.
\end{equation}
Hence there exists some $T_0>0$ such that for all $(x,t)\in\mathbb{R}\times[T_0,+\infty)$, it holds that
$$
-\frac{\t_-}{2}<\zeta<\frac{\t_+}{2}.
$$
This directly yields, for all $(x,t)\in\mathbb{R}\times[T_0,+\infty)$,
\begin{equation}\label{theta-}
\di \theta=\Theta+\zeta>\theta_- -\frac{\theta_-}{2}=\frac{\theta_-}{2},
\end{equation}
and
\begin{equation}
\di \theta=\Theta+\zeta<\theta_+ +\frac{\theta_+}{2}=\frac{3\theta_+}{2}.
\end{equation}

Finally, it follows from $\eqref{ns}_3$ that,
\begin{equation}
\di \t_t+\frac{R\t}{c_{\nu}}\frac{u_x}{v}-\frac{\k}{c_{\nu}}\frac{\t_{xx}}{v}+\frac{\k}{c_{\nu}}\frac{\t_xv_x}{v^2}\geq 0.
\end{equation}
Define
\begin{equation}
\di \bar\t=\t \exp\left(\frac{R}{c_{\nu}}\int_0^t\left\|\frac{u_x}{v}\right\|_{L^{\infty}}ds\right).
\end{equation}
We find
\begin{equation}
\di \bar\t_t-\frac{\k}{c_{\nu}}\frac{\bar\t_{xx}}{v}+\frac{\k}{c_{\nu}}\frac{\bar\t_xv_x}{v^2}\geq
\frac{R\bar\t}{c_{\nu}}\Big(\left\|\frac{u_x}{v}\right\|-\frac{u_x}{v}
\Big)\geq 0.
\end{equation}
By the minimum principle of the parabolic equation, we obtain
\begin{equation}
\di\inf_{x,t}\bar\t\geq\inf_{x,t}\bar\t\Big|_{t=0}=\inf_{x\in\mathbb{R}}\t_0\geq m^{-1}_0,
\end{equation}
which directly yields
\begin{equation}
\begin{array}{ll}\label{thata-lower}
\di\inf_{x,t}\t(x,t)&\di\geq m^{-1}_0e^{-\frac{R}{c_{\nu}}\int_0^t\|\frac{u_x}{v}\|_{L^{\infty}}ds}\\
&\di\geq m_0^{-1}e^{-C_0\int_0^t\|\psi_x\|_{L^{\infty}}+\|U_x\|_{L^{\infty}}ds}\\
&\di\geq m_0^{-1} e^{-C_0\int_0^t\|\psi_x\|^{1/2}\|\psi_{xx}\|^{1/2}}e^{-C_0\int_0^t\|U_x\|_{L^{\infty}}ds}\\
&\di\geq m_0^{-1} e^{-C_0\int_0^t\|\psi_{xx}\|^{1/2}}e^{-C_0\d t}\\
&\di\geq m_0^{-1} e^{-C_0(t^{3/4}+\d t)}\geq C_0e^{-C_0t}.
\end{array}
\end{equation}
Thus from \eqref{theta-}, we have $\t>\min\{\frac{\t_-}{2},C_0e^{-C_0T_0}\}$
for all $(x,t)\in\mathbb{R}\times(0,+\infty)$.  By Lemma \ref{basic-lemma},
\ref{v-bound}-\ref{high-order}, Proposition \ref{prop} is completed.
\hfill $\Box$

\


\section{Proof of Theorem \ref{theorem2}}
\setcounter{equation}{0}

It is sufficient to show the same a priori estimate as Proposition \ref{prop}. Noticing that
$(V_{\pm}^r,U_{\pm}^r,\T_{\pm}^r)$ satisfies Euler system \eqref{euler} and $(V^{cd},U^{cd},\T^{cd})$
satisfies $\eqref{ns}_1$ \eqref{Theta t}, we rewrite the Cauchy problem \eqref{ns}\eqref{initial} as
\begin{equation}\label{perturb2}
\left\{
\begin{array}{ll}
\di \phi_t-\psi_x=0,\\
\di \psi_t+\left(p-P\right)_x=\mu\left(\frac{u_x}{v}-\frac{U_x}{V}\right)_x+F,\\
\di c_{\nu}\zeta_t+pu_x-PU_x=\kappa\left(\frac{\t_x}{v}-\frac{\Theta_x}{V}\right)_x
+\mu\left(\frac{u_x^2}{v}-\frac{U^2_x}{V}\right)+G,\\
\di(\phi,\psi,\zeta)(\pm\infty,t)=0,\\
\di (\phi,\psi,\zeta)(x,0)=(\phi_0,\psi_0,\zeta_0)(x),\quad x\in \mathbb{R},
\end{array}
\right.
\end{equation}
where
$$
P=\frac{R\T}{V},\quad P_{\pm}=\frac{R\T_{\pm}^r}{V_{\pm}^r},
$$
$$
F=(P_-+P_+-P)_x+\left(\mu\frac{U_x}{V}\right)_x-U_t^{cd},
$$
and
$$
\begin{array}{ll}
G=(p^m-P)U^{cd}_x+(P_--P)(U_-^r)_x+(P_+-P)(U_+^r)_x\\
\di\quad\quad+\mu\frac{U_x^2}{V}+\kappa\left(\frac{\T_x}{V}-\frac{\T_x^{cd}}{V^{cd}}\right)_x.
\end{array}
$$

Similar to Lemma \ref{basic-lemma}, the following key estimate holds.
\begin{lemma}\label{basic2}
For $(\phi,\psi,\zeta)\in X([0,T])$, we assume \eqref{same-order} holds, then there exist some positive
constants $C_0$ and $\d_0$ such that if $\d<\d_0$, it follows that for $t\in[0,T]$,
\begin{equation}
\begin{array}{ll}
\di \Big\|\left(\psi,\sqrt{\Phi\left(\frac{v}{V}\right)},\sqrt{\Phi\left(\frac{\t}{\Theta}\right)}\right)(t)\Big\|^2
+\int_0^t\int\left(\frac{\psi_x^2}{\t v}+\frac{\zeta_x^2}{\t^2 v}\right)dxds\\[3mm]
\di+\int_0^t\int P\left(\Phi\left(\frac{\t V}{\T v}\right)+\g\Phi\left(\frac{v}{V}\right)\right)\big((U_-^r)_x+(U_+^r)_x\big)dxds
\di\leq C_0,
\end{array}
\end{equation}
where $C_0$ denotes a constant depending only on $\mu$, $\kappa$, $R$, $c_{\nu}$,
$v_{\pm}$, $u_{\pm}$, $\t_{\pm}$  and $m_0$.
\end{lemma}

\textbf{Proof}: First, multiplying $\eqref{perturb2}_2$ by $\psi$ leads to
\begin{equation}\label{4.17}
\begin{array}{ll}
\di \left(\frac{\psi^2}{2}\right)_t+\left[(p-P)\psi-\mu\left(\frac{u_x}{v}-\frac{U_x}{V}\right)\psi\right]_x
-\frac{R\zeta}{v}\psi_x\\
\di-R\T\left(\frac{1}{v}-\frac{1}{V}\right)\phi_t+\mu\frac{\psi_x^2}{v}
+\mu\left(\frac{1}{v}-\frac{1}{V}\right)U_x\psi_x=F\psi.
\end{array}
\end{equation}

Next, we multiply $\eqref{perturb2}_3$ by $\z\t^{-1}$ to get
\begin{equation}\label{4.18}
\begin{array}{ll}
\di\frac{R}{\g-1}\frac{\z\z_t}{\t}-\left[\k\left(\frac{\t_x}{v}-\frac{\T_x}{V}\right)\frac{\z}{\t}\right]_x
+\frac{R\z}{v}\psi_x+(p-P)U_x\frac{\z}{\t}\\[2mm]
\di+\k\frac{\T\z_x^2}{\t^2v}-\k\frac{\z_x\z\T_x}{\t^2v}-\k\frac{\phi\T_x\T\z_x}{\t^2vV}
+\k\frac{\phi\z\T_x^2}{\t^2vV}\\[2mm]
\di-\mu\frac{\z\psi_x^2}{\t v}-2\mu\frac{\psi_xU_x\z}{\t v}+\mu\frac{\phi\z U_x^2}{\t vV}=G\frac{\z}{\t}.
\end{array}
\end{equation}
Noticing that
\begin{equation}\label{4.19}
\di-R\T\left(\frac{1}{v}-\frac{1}{V}\right)\phi_t=\left[R\T\Phi\left(\frac{v}{V}\right)\right]_t
-R\T_t\Phi\left(\frac{v}{V}\right)+\frac{P\phi^2}{vV}V_t,
\end{equation}
\begin{equation}
\di\frac{\z\z_t}{\t}=\left[\T\Phi\left(\frac{\t}{\T}\right)\right]_t+\T_t\Phi\left(\frac{\T}{\t}\right),
\end{equation}
\begin{equation}
\begin{array}{ll}
-R\T_t&\di=(\g-1)P_-(U_-^r)_x+(\g-1)P_+(U_+^r)_x-p^mU_x^{cd}\\[2mm]
&\di=(\g-1)P(U_-^r)_x+(\g-1)P(U_+^r)_x+(\g-1)(P_--P)(U_-^r)_x\\[2mm]
&\di+(\g-1)(P_+-P)(U_+^r)_x-p^mU_x^{cd},
\end{array}
\end{equation}
and
\begin{equation}
\begin{array}{ll}
\di-R\T_t\Phi\left(\frac{v}{V}\right)+\frac{P\phi^2}{vV}V_t+\frac{R}{\g-1}\T_t\Phi\left(\frac{\T}{\t}\right)
+(p-P)U_x\frac{\z}{\t}\\
=Q_1\big((U_-^r)_x+(U_+^r)_x\big)+Q_2,
\end{array}
\end{equation}
where
\begin{equation}
\begin{array}{ll}
Q_1&\di=(\g-1)P\Phi\left(\frac{v}{V}\right)+\frac{P\phi^2}{vV}-P\Phi\left(\frac{\T}{\t}\right)
+\frac{\z}{\t}(p-P)\\
&\di=P\left(\Phi\left(\frac{\t V}{\T v}\right)+\g\Phi\left(\frac{v}{V}\right)\right),
\end{array}
\end{equation}
and
\begin{equation}
\begin{array}{ll}\label{q2}
Q_2&\di=U_x^{cd}\left(\frac{P\phi^2}{vV}-p^m\Phi\left(\frac{v}{V}\right)+\frac{p^m}{\g-1}\Phi\left(\frac{\T}{\t}\right)
+\frac{\z}{\t}(p-P)\right)\\[3mm]
&\di+(\g-1)(P_--P)(U_-^r)_x\left(\Phi\left(\frac{v}{V}\right)-\frac{1}{\g-1}\Phi\left(\frac{\T}{\t}\right)\right)\\[3mm]
&\di+(\g-1)(P_+-P)(U_+^r)_x\left(\Phi\left(\frac{v}{V}\right)-\frac{1}{\g-1}\Phi\left(\frac{\T}{\t}\right)\right).
\end{array}
\end{equation}
Combining \eqref{4.18} and \eqref{4.19}, it follows from \eqref{4.19}-\eqref{q2} that
\begin{equation}\label{2-basic}
\begin{array}{ll}
\di\left(\frac{\psi^2}{2}+R\Theta\Phi\left(\frac{v}{V}\right)+\frac{R}{\g-1}\Theta\Phi\left(\frac{\theta}{\Theta}\right)\right)_t
+\frac{\mu\Theta}{\theta v}\psi_x^2+\frac{\kappa\Theta}{\theta^2 v}\zeta_x^2\\
\di+H_{x}+Q_1\big((U_-^r)_x+(U_+^r)_x\big)+\tilde Q=F\psi+G\frac{\zeta}{\theta}
\end{array}
\end{equation}
with $H$ the same as in \eqref{H}, and
\begin{equation}\label{tilde-q}
\begin{array}{ll}
\tilde Q&\di=Q_2-\k\frac{\z_x\z\T_x}{\t^2v}-\k\frac{\phi\T_x\T\z_x}{\t^2vV}
+\k\frac{\phi\z\T_x^2}{\t^2vV}\\[2mm]
&\di-\mu\frac{\phi U_x\psi_x}{vV}-2\mu\frac{\psi_xU_x\z}{\t v}+\mu\frac{\phi\z U_x^2}{\t vV}.
\end{array}
\end{equation}
Recalling (iii) in Lemma \ref{rare-pro}, we can compute
\begin{equation}
\begin{array}{ll}
\di|(P_--P)(U_-^r)_x|\\
\di\leq C\big(|\T^{cd}-\t_-^m|+|\T^{r}_+-\t_+^m|+|V^{cd}-v_-^m|+|V^{r}_+-v_+^m|\big)|(U_-^r)_x|\\
\di\leq C\big(|\T^{cd}-\t_-^m|+|\T^{r}_+-\t_+^m|+|V^{cd}-v_-^m|+|V^{r}_+-v_+^m|\big)\big|_{\Omega_-}
+C|(U_-^r)_x|\big|_{\Omega_c\cap\Omega_+}\\
\di\leq C\d e^{-c_0(|x|+t)},
\end{array}
\end{equation}
which leads to
\begin{equation}
\di|Q_2|\leq C(M)|U_x^{cd}|(\phi^2+\zeta^2)+C(M)\d e^{-c_0(|x|+t)}(\phi^2+\zeta^2)
\end{equation}
and
\begin{equation}
\di|\tilde Q|\leq |Q_2|+\frac{\mu\T\psi_x^2}{4\t v}+\frac{\k\T\zeta_x^2}{4\t^2 v}
+C(M)(\phi^2+\zeta^2)(\T_x^2+U_x^2).
\end{equation}
Note that
\begin{equation}
\begin{array}{ll}
(\phi^2+\zeta^2)(\T_x^2+U_x^2)\leq C(\phi^2+\zeta^2)((\T_x^{cd})^2+(U^r_-)_x^2+(U^r_+)_x^2)\\
\quad \leq C(1+t)^{-1}(\phi^2+\zeta^2)\e^{-\frac{c_1x^2}{1+t}}+C\d Q_1((U^r_-)_x+(U^r_+)_x).
\end{array}
\end{equation}
Following the same calculations as in \cite{Huang-Li-matsumura}, it holds that
\begin{equation}
\di \|(F,G)\|_{L^1}\leq C\d^{1/8}(1+t)^{-7/8}.
\end{equation}
Then we have
\begin{equation}\label{int-F}
\begin{array}{ll}
\di\int_0^t\int\left(F\psi+G\frac{\z}{\t}\right)dxds\leq C(M)\int_0^t\|(F,G)\|_{L^1}\|(\psi,\zeta)\|_{L^{\infty}}ds\\
\di\quad\leq C(M)\d^{1/8}\int_0^t(1+s)^{-7/8}\|(\psi,\zeta)\|^{1/2}\|(\psi_x,\zeta_x)\|^{1/2}ds\\
\di\quad\leq \int_0^t\int\left(\frac{\mu\T\psi_x^2}{4\t v}+\frac{\k\T\zeta_x^2}{4\t^2v}\right)dxds\\
\di\quad+C(M)\d^{1/6}\int_0^t(1+s)^{-7/6}\left(1+\left\|\left(\psi,\sqrt{\Phi\left(\frac{\t}{\T}\right)}\right)\right\|^2\right)ds.
\end{array}
\end{equation}
Integrating \eqref{2-basic} over $\mathbb{R}\times(0,t)$ and using Gronwall's inequality,
we deduce from \eqref{tilde-q}-\eqref{int-F} that
\begin{equation}\label{7.18}
\begin{array}{ll}
\di \Big\|\left(\psi,\sqrt{\Phi\left(\frac{v}{V}\right)},\sqrt{\Phi\left(\frac{\t}{\Theta}\right)}\right)(t)\Big\|^2
+\int_0^t\int\left(\frac{\psi_x^2}{\t v}+\frac{\zeta_x^2}{\t^2 v}\right)dxds\\[3mm]
\di+\int_0^t\int Q_1\big((U_-^r)_x+(U_+^r)_x\big)dxds \leq C_0+C(M)\d^{1/6}\\[3mm]
\di+C(M)\d\int_0^t(1+s)^{-1}\int_{\mathbb{R}}(\phi^2+\zeta^2)e^{-\frac{c_1x^2}{1+s}}dxds.
\end{array}
\end{equation}
Finally, due to the fact that
\begin{equation}\label{7.19}
\begin{array}{ll}
\di\int_0^t\int(\phi^2+\psi^2+\zeta^2)w^2dxds
&\di\leq C(M)+C(M)\int_0^t\|(\phi_x,\psi_x,\zeta_x)\|^2ds\\
&\di+C(M)\int_0^t\int(\phi^2+\zeta^2)\big((U_-^r)_x+(U_+^r)_x\big)dxds\\
&\di\leq C(M)+C(M)\int_0^t\int\left(\frac{\t\phi_x^2}{v^3}+\frac{\psi_x^2}{\t v}
+\frac{\zeta_x^2}{\t^2v}\right)dxds\\
&\di+C(M)\int_0^t\int Q_1\big((U_-^r)_x+(U_+^r)_x\big)dxds,
\end{array}
\end{equation}
whose proof can be found in \cite{Huang-Li-matsumura} and $w$ is defined as in Lemma \ref{hlm}, substituting \eqref{7.19} into \eqref{7.18} and choosing $\d$ suitable small imply
\begin{equation}
\begin{array}{ll}
\di \Big\|\left(\psi,\sqrt{\Phi\left(\frac{v}{V}\right)},\sqrt{\Phi\left(\frac{\t}{\Theta}\right)}\right)(t)\Big\|^2
+\int_0^t\int\left(\frac{\psi_x^2}{\t v}+\frac{\zeta_x^2}{\t^2 v}\right)dxds\\[3mm]
\di+\int_0^t\int Q_1\big((U_-^r)_x+(U_+^r)_x\big)dxds
\leq C_0+C(M)\d\int_0^t\int\frac{\t\phi_x^2}{v^3}dxds.
\end{array}
\end{equation}
Thus we can finish the proof of Lemma \ref{basic2} in the similar way as in Lemma \ref{basic-lemma}. We omit the details for brevity.
\hfill $\Box$

It is easy to check that the other estimates for single viscous contact wave still hold for the case in which the composite waves are the combination of viscous contact wave with rarefaction waves.
Thus, we complete the proof of Proposition \ref{prop}, and finally prove Theorem \ref{theorem2}.




\small
\begin{thebibliography}{99}

\bibitem{goodmann}
\newblock J. Goodman,
\newblock\emph{ Nonlinear asymptotic stability of viscous shock profiles for conservation laws},
\newblock Arch. Ration. Mech. Anal., \textbf{95} (1986), 325--344.

\bibitem{xiao-liu}
\newblock L. Hsiao, T. Liu,
\newblock\emph{Nonlinear diffusive phenomena of nonlinear hyperbolic systems},
\newblock Chinese Ann. Math. Ser. B, \textbf{14} (1993), 465--480.

\bibitem{Hong}
\newblock H. H. Hong,
\newblock\emph{Global stability of viscous contact wave for 1-D compressible Navier-Stokes equations},
\newblock J. Differential Equations, \textbf{252} (2012), 3482--3505.


\bibitem{Huang-Li-matsumura}
\newblock F. M. Huang, J. Li and  A. Matsumura,
\newblock\emph{ Asymptotic stability of combination of viscous contact wave with rarefaction waves for one-dimensional compressible Navier-Stokes system},
\newblock  Arch. Ration. Mech. Anal.,\textbf{197} (2010), 89--116.

\bibitem{Huang-matsumura}
\newblock F. M. Huang,   A. Matsumura,
\newblock\emph{ Stability of a composite wave of two viscous shock waves for the full compressible Navier-Stokes equation},
\newblock  Comm. Math. Phys., \textbf{289} (2009), 841--861.


\bibitem{huang-ma-shi}
\newblock F. H. Huang, A. Matsumura and X. D. Shi,
\newblock\emph{On the stability of contact discontinuity for compressible Navier-Stokes equations with free boundary},
\newblock Osaka J. Math., \textbf{41} (2004), 193--210.



\bibitem{Huang-Matsumura-Xin}
\newblock F. M. Huang, A. Matsumura and Z. P. Xin,
\newblock \emph{Stability of contact discontinuities for the 1-D compressible
Navier-Stokes equations},
\newblock  Arch. Ration. Mech. Anal., \textbf{179} (2006), 55--77.





\bibitem{Huang-Xin-Yang}
\newblock F. M. Huang, Z. P. Xin and T. Yang,
\newblock \emph{Contact discontinuities with general perturbation for gas motion},
\newblock Adv. Math., \textbf{219} (2008), 1246--1297.

\bibitem{Huang-Yang}
\newblock F. M. Huang, T. Yang,
\newblock\emph{Stability of contact discontinuity for the Boltzmann equation},
\newblock J. Differential Equations, \textbf{229} (2006), 698--742.

\bibitem{huang-zhao}
\newblock F. M. Huang, H. J. Zhao,
\newblock\emph{On the global stability of contact discontinuity for compressible
Navier-Stokes equations},
\newblock Rend. Sem. Mat. Univ. Padova., \textbf{109} (2003), 283--305.

%

\bibitem{jiang}
\newblock S. Jiang,
\newblock \emph{Large-time behavior of solutions to the equations of a one-dimensional
viscous polytropic ideal gas in unbounded domains},
\newblock Comm. Math. Phys., \textbf{200} (1999), 181--193.

\bibitem{jiang-2}
\newblock S. Jiang,
\newblock \emph{Remarks on the asymptotic behaviour of solutions to the compressible
Navier-Stokes equations in the half-line},
\newblock Proc. Roy. Soc. Edinburgh Sect. A,   \textbf{132} (2002), 627--638.


\bibitem{KM}
\newblock S. Kawashima ,  A. Matsumura,
\newblock \emph{Asymptotic stability of traveling wave solutions of systems for one-dimensional gas
motion},
\newblock Comm. Math. Phys., \textbf{101}  (1985) 97--127.

\bibitem{kazhikhov}
\newblock A. V. Kazhikhov,
\newblock \emph{Cauchy problem for viscous gas equations},
\newblock Siberian Math. J., \textbf{23} (1982), 44--49.

\bibitem{kazhi-she}
\newblock A. V. Kazhikhov, V. V. Shelukin,
\newblock\emph{Unique global solution with respect to time
of initial boundary value problems for one-dimensional equations of a viscous gas},
\newblock J. Appl. Math. Mech. \textbf{41}  (1977), 273--282.

\bibitem{Lax}P. Lax,  \textit{Hyperbolic systems of conservation laws, II,} Comm. Pure Appl. Math.,
{\bf 10} (1957), 537-566.

\bibitem{L-L}
J. Li, Z. L. Liang,
\textit{Some Uniform Estimates and Large-Time Behavior for
One-Dimensional Compressible Navier-Stokes System in
Unbounded Domains with Large Data},  http://arxiv.org/abs/1404.2214.


\bibitem{Liu-1985} T. P. Liu, \textit{Nonlinear stability of shock waves for
viscous conservation laws}, Mem. Amer. Math. Soc., \textbf {56}
(1985), 1--108.

\bibitem{Liu-Xin}
\newblock T. P. Liu,  Z. P.  Xin,
\newblock \emph{Pointwise decay to contact discontinuities for systems of viscous conservation laws},
\newblock  Asian J. Math., \textbf{1} (1997) 34--84.


\bibitem{MN-85}
\newblock A. Matsumura, K. Nishihara,
\newblock \emph{On the stability of traveling wave solutions of a one-dimensional model system for
compressible viscous gas},
\newblock  Japan J. Appl. Math., \textbf{2}  (1985) 17--25.

\bibitem{MN-86}
\newblock A. Matsumura, K. Nishihara,
\newblock \emph{Asymptotics toward the rarefaction wave of the solutions of
a one-dimensional model system for compressible viscous gas},
\newblock  Japan J. Appl. Math., \textbf{3} (1986), 1--13.

\bibitem{MN-92}
\newblock A. Matsumura, K. Nishihara,
\newblock \emph{Global stability of the rarefaction wave of a onedimensional
model system for compressible viscous gas},
\newblock  Comm. Math. Phys., \textbf{144} (1992), 325--335.

\bibitem{nishi-yang-zhao}
K. Nishihara, T. Yang, H. J. Zhao,
\emph{Nonlinear stability of strong rarefaction waves for compressible Navier-Stokes equations},
SIAM J. Math. Anal., \textbf{35}  (2004),  1561--1597.

\bibitem{Smoller} 
     \newblock J. Smoller,
     \newblock ``Shock Waves and Reaction-Diffusion Equations," New
     York: Springer, 1994.

\bibitem{szepessy-xin}
A. Szepessy and Z. P. Xin,
\emph{Nonlinear stability of viscous shock waves},
Arch. Ration. Mech. Anal., \textbf{122} (1993), 53--103.


\bibitem{Xin2}
\newblock Z. P. Xin,
\newblock \emph{On nonlinear stability of contact discontinuities. Hyperbolic problems: theory, numerics, applications  Stony Brook, NY, 1994},
\newblock  World Sci. Publishing, River Edge, NJ, (1996),  249--257.





\end{thebibliography}
\end{document}